\documentclass[11pt]{article}
\usepackage{amssymb}
\usepackage{amssymb}
\usepackage{amssymb}
\usepackage{amsfonts}
\usepackage{float}
\usepackage{amsfonts}
\def\oneandahalf{\par\baselineskip=18pt}
\def\standard{\oddsidemargin=0in
              \evensidemargin=0in
              \topmargin =-.5in
              \textheight=9.0in
              \textwidth=6.5in}

{\end{list}}
\newenvironment{hangref}{\begin{list}{}{\setlength{\itemsep}{0pt}%
\setlength{\parsep}{0pt}\setlength{\leftmargin}{+\parindent}%
\setlength{\itemindent}{-\parindent}}}{\end{list}}
\def\E{{\rm E}}
\def\Var{{\rm Var}}

\def\bG{{\bf G}}

\def\bt{{\bf t}}

\def\bX{{\bf X}}
\def\bx{{\bf x}}
\def\bY{{\bf Y}}
\def\by{{\bf y}}
\def\bZ{{\bf Z}}
\def\bz{{\bf z}}
\def\b1{{\bf 1}}

\def\blot{\quad {$\vcenter{\vbox{\hrule height.4pt
             \hbox{\vrule width.4pt height.9ex \kern.9ex \vrule
width.4pt}
             \hrule height.4pt}}$}}
\usepackage{amssymb,latexsym}
\usepackage{amsmath}
\usepackage{graphicx}
\usepackage{amsmath,amsfonts,amssymb,amsbsy,amsthm}
\usepackage{enumerate}
\standard
\usepackage{float}
\usepackage{color}
\newtheorem{theorem}{Theorem}
\newtheorem{lemma}{Lemma}

\newtheorem{pro}{Proposition}
\newtheorem{remark}{Remark}
\newtheorem{assumption}{Assumption}
\newtheorem{condition}{Condition}
\newtheorem{definition}{Definition}
\newtheorem{exmp}{Example}
\begin{document}

\normalsize
\title{Conditional Monte Carlo: A Change-of-Variables Approach}

\author{
        Guiyun Feng                         \\
        {\footnotesize
        Department of Industrial \& Systems Engineering, University of
        Minnesota}\\ \footnotesize{111 Church Street S.E., Minneapolis, MN 55455, USA}\\
        \\
        Guangwu Liu                            \\
        {\footnotesize
        Department of Management Sciences, City University of Hong Kong}\\
        \footnotesize{Tat Chee Avenue, Kowloon, Hong Kong} }

\date{\footnotesize{\today}}

\maketitle

\begin{abstract}

Conditional Monte Carlo (CMC) has been widely used for sensitivity
estimation with discontinuous integrands as a standard simulation
technique. A major limitation of using CMC in this context is that
finding conditioning variables to ensure continuity and tractability
of the resulting conditional expectation is often problem dependent
and may be difficult. In this paper, we attempt to circumvent this
difficulty by proposing a change-of-variables approach to CMC,
leading to efficient sensitivity estimators under mild conditions
that are satisfied by a wide class of discontinuous integrands.
These estimators do not rely on the structure of the simulation
models and are less problem dependent. The value of the proposed
approach is exemplified through applications in sensitivity
estimation for financial options and gradient estimation of
chance-constrained optimization problems.

\end{abstract}

\oneandahalf

\section{Introduction}

Conditional Monte Carlo (CMC) is a standard simulation technique
that has been widely discussed by many simulation textbooks; see,
e.g., Law and Kelton (2000) and Asmussen and Glynn (2007). When
estimating the expectation of a random performance, the basic idea
of CMC is to use conditional expectation of the performance, rather
than the random performance itself, as an estimator, where
appropriate conditioning variables are chosen to obtain the
conditional expectation. It is well known that the CMC estimator has
a smaller variance, as guaranteed by the law of total variance.

In this paper, our main interest is focused on the use of CMC in
sensitivity estimation of an expectation, another context in which
CMC plays an important role. In particular, we consider a setting
that the integrand of the expectation is discontinuous, under which
sensitivity estimation is challenging, and has received a
significant amount of attention in the simulation community in
recent years. In this setting, CMC works in conjunction with
infinitesimal perturbation analysis (IPA, also known as the pathwise
method) that suggests interchanging the order of differentiation and
expectation. Typically, IPA does not work when the integrand is
discontinuous, because the interchange is not valid due to
discontinuity. To fix this issue, CMC suggests finding appropriate
conditioning variables and then applying IPA on the conditional
expectation that is often continuous. The intuition behind is that
taking conditional expectation can often ``integrate out"
discontinuity and smooth the integrand. This idea was first proposed
by Gong and Ho (1987) and Suri and Zazanis (1988), and is also
referred to as smoothed perturbation analysis; see Fu and Hu (1997)
for a monograph on detailed treatments, and Wang at al. (2009) and
Fu et al. (2009) for some of its recent applications.

A major limitation of CMC for sensitivity estimation is how to find
conditioning variables such that the conditional expectation is
smooth and easily computable. Finding such conditioning variables is
problem dependent and could be difficult for some cases. For
instance, when discontinuity comes from an indicator function that
equals $1$ if the maximum of a random vector is smaller than a
threshold and $0$ otherwise, it may be difficult to find
conditioning variables to ensure continuity of the resulting
conditional expectation. This limitation motivates us to reexamine
CMC from other perspectives, aiming to circumvent the difficulty on
finding conditioning variables.

In this paper, we study CMC from a change-of-variables perspective,
and we call it {\it a change-of-variables approach}. Intuitively, it
proceeds by constructing a one-to-one mapping and applying a change
of variables, followed by taking iterated integrations. The approach
is appealing when the discontinuous integrand involves an indicator
function. In this setting, we want to estimate the sensitivities of
the expectation of
\begin{equation}\label{eqn:IndForm}
l(\bX)\cdot 1_{\left\{h(\bX)\le \xi\right\}},
\end{equation}
for continuous functions $l(\cdot)$ and $h(\cdot)$, where $\bX$ is a
random vector that captures the randomness of the simulation, and
$\xi$ is a given constant. A one-to-one mapping $\bX\mapsto
(u(\bX),h(\bX))$ can be constructed with a vector function $u$, and
a change-of-variables formula represents the expectation of
(\ref{eqn:IndForm}) as a double integral taken over the supports of
$u(\bX)$ and $h(\bX)$, respectively. Then integrating along the
dimension of $h(\bX)$ over $(-\infty,\xi)$ removes the indicator and
produces a smooth integrand, which enables the use of IPA for
sensitivity estimation.

Though the intuition of the change-of-variables approach is
straightforward, its theoretical justification under a general
setting is not trivial, because the domain and image sets of the
mapping may  have different Euclidean dimensions. For instance,
when $\bX\triangleq (X_1,\dots,X_m)$ and
$h(\bX)=\max(X_1,\dots,X_m)$, a useful one-to-one mapping that we
may construct is $\bX\mapsto (\bX/h(\bX),h(\bX))$. While the domain
 of this mapping is a subset of $\mathbb{R}^m$, its image  is
a subset of $\mathbb{R}^{m+1}$, which is indeed an $m$-dimensional
manifold in $\mathbb{R}^{m+1}$. Change-of-variables formulas for
such cases involve an extension of Lebesgue measure, and more
generally, geometric measure theory; see Federer (1996) for a
monograph. In this paper, we provide theoretical underpinnings of
the approach under a general setting, and discuss how it can be applied to develop sensitivity estimators. %Provided that the $h(\cdot)$ function satisfies certain
%smoothness conditions, we show that the approach leads to new and
%efficient estimators which do not rely on the structure of the
%simulation model and are thus less problem dependent.

Sensitivity estimation finds applications in a wide range of areas
in operations research. For instance, it can be used in simulation
optimization to estimate gradients that serve as key inputs to many
gradient-based optimization algorithms. In financial applications,
it estimates hedging parameters such as delta and gamma, which play
important roles in risk management of financial securities. There
has been a vast literature on sensitivity estimation, and various
methods have been proposed, traditional ones including
finite-difference approximations, IPA (see, e.g., Ho and Cao 1983
and Broadie and Glasserman 1996), the likelihood ratio method (see,
e.g., Glynn 1987 and L'Ecuyer 1990), the weak derivative method
(see, e.g., Pflug 1988 and Pflug and Weisshaupt 2005), and the
Malliavin calculus method (see, e.g., Bernis et al. 2003 and Chen
and Glasserman 2007). In recent years, sensitivity estimation for
expectations with discontinuous integrands has received a
significant amount of attention among simulation researchers and
several methods have been proposed. Lyuu and Teng (2011) showed that
the sensitivity can be written as an integral taken over an
appropriate subset, and suggested using importance sampling to
estimate the integral. Liu and Hong (2011) showed that the
sensitivity is a summation of two terms with the latter one
involving a conditional expectation and a density, and proposed a
kernel smoothing method to estimate the second term. Wang et al.
(2012) proposed the so-called SLRIPA method that moves the parameter
of interest out of the indicator function to smooth the integrand
and enables the use of IPA. It unifies the likelihood ratio method
and the IPA method in certain sense. Chan and Joshi (2013) suggested
a novel bumping on sample paths in a way that discontinuous points
are eliminated. Essentially their approach relies on an appropriate
change of variables on the sample paths. Along the line of Liu and
Hong (2011), Tong and Liu (2016) proposed an importance sampling
method to estimate conditional expectations, leading to unbiased
estimators of the sensitivities.

%One of
%the simulation methods that deal with the sensitivity
%estimation problem builds upon the interchange of the
%differentiation and expectation, which is usually referred to as the
%infinitesimal perturbation analysis (IPA) or pathwise method. For
%the above problem, however, such an interchange is not valid because
%the performance is not continuous in $\theta$, while continuity is a
%crucial requirement for IPA. To fix this issue, CMC suggests
%choosing a random vector $\bY$ and applying IPA on the conditional
%expectation $\E\left[\left.l(\bX)\cdot1_{\{h(\bX)\le
%\xi\}}\right|\bY\right]$, based on an intuition that conditioning
%can often ``integrate out" the discontinuity and smooth the
%integrand. This idea was first proposed by Gong and Ho (1987) and
%Suri and Zazanis (1988) and is also referred to as smoothed
%perturbation analysis; see also Fu and Hu (1997) for a monograph on
%detailed treatments of CMC and various applications. However, how to
%choose $\bY$ such that the conditional expectation
%$\E\left[\left.l(\bX)\cdot1_{\{h(\bX)\le \xi\}}\right|\bY\right]$ is
%continuous and easily computable is problem dependent. As far as we
%are aware, it is clear how to choose appropriate $\bY$ such that the
%conditional expectation is continuous in $\theta$ and computable for
%some important cases, e.g., when $h(\bX)=\max(X_1,\ldots,X_m)$ as in
%sensitivity estimation of the prices of many barrier and lookback
%financial options. This gap motivates us to reexamine CMC from other
%perspectives.

As a remark, we would like to point out that the change-of-variables
idea is not new for Monte Carlo simulation. Dating back to 1950s, a
change-of-variables argument was discussed by Hammersley (1956) in
the context of computing a conditional expectation at a fixed point,
and further elaborated by Wendel (1957) from a group-theoretic
aspect. This line of research, however, received little attention
afterwards and its theoretical underpinnings were underdeveloped. %In
%this paper, we study in depth the change-of-variables approach in a
%general CMC setting, and provide theoretical underpinnings, based on
%which we develop new and efficient estimators for sensitivities of
%expectations with discontinuous integrands.

To summarize, we make the following contributions in this paper.
\begin{itemize}
\item We develop a change-of-variables framework of CMC for general integrands and provide theoretical underpinnings.
This framework is adapted to Hausdorff measure that offers great
flexibility in construction of CMC estimators.
\item We propose sensitivity estimators for expectations with discontinuous integrands of the form
(\ref{eqn:IndForm}). These estimators require only the probability
density function of $\bX$ as an input and certain smoothness
conditions on the function $h$. They are less problem dependent, and
work for a wide class of $h$'s.
\item We study two applications, including sensitivity estimation for financial options with discontinuous payoffs and
gradient estimation of chance-constrained optimization problems, to
illustrate how the proposed approach may lead to new and efficient
estimators.
\end{itemize}

The remainder of the paper is organized as follows. We formulate the
problem in Section \ref{sec:Problem}, and introduce the main idea of
the change-of-variables approach in Section \ref{sec:Perspevtive}.
Section \ref{sec:Sensitivity} shows how to use the
change-of-variables approach for sensitivity estimation with
discontinuous integrands. Two applications with numerical examples
are considered in Section \ref{sec:Application} to exemplify the
value of the proposed approach, followed by concluding remarks in
Section \ref{sec:Conclusions}. Proofs of the main results are
provided in the appendix, while lengthy details of the estimators
used in numerical examples are put in an online supplement.

\vspace{10pt}

\noindent\textbf{Notation}. Throughout the paper, we let boldface
and italic letters denote vectors and vector elements respectively,
and upper-case letters denote random counterparts of the lower-case
ones. For a mapping $u$, $u\{A\}$ denotes the set of images of a set
$A$ in the domain, while $u(\bx)$ denotes the image of an element
$\bx\in A$.

\section{Problem Formulation}\label{sec:Problem}

Consider a sensitivity estimation problem for an expectation with a
discontinuous integrand of the form
\begin{equation*}
l(\bX)\cdot1_{\left\{h(\bX)\le \xi\right\}},
\end{equation*}
where $l$ and $h$ are continuous functions, $\xi$ is a given
constant, and $\bX=(X_1,\dots,X_m)$ is a random vector that captures
the randomness of the simulation. This form of integrands has
received a significant amount of attention in the simulation society
in recent years; see, e.g., Lyuu and Teng (2011), Hong and Liu
(2011), Wang et al. (2012), Chan and Joshi (2013), and Tong and Liu
(2016).

Throughout the paper, we assume that $\bX$ has a known density
function.
\begin{assumption}\label{assum:Density}
The random vector $\bX$ has a density function $f(\bx)$ on a support
$\Omega\subset \mathbb{R}^m$.
\end{assumption}
This assumption can be further relaxed to a requirement of a
conditional density of $\bX$ given a random vector $\bG$; see
Section \ref{sec:VG} of the online supplement for such an example.

Suppose the sensitivity of interest is taken with respect to
(w.r.t.) $\theta$, a parameter on which $\bX$ may depend. Without
loss of generality we assume that $\theta$ is one-dimensional and
$\theta\in\Theta$, where $\Theta$ is an open set. If $\theta$ is
multidimensional, one may treat each dimension as a one-dimensional
parameter while fixing other dimensions as constants.

Explicitly accounting for its dependence on $\theta$, the quantity
we want to estimate is written as
\begin{equation}\label{eqn:DiscontTheta}
\gamma'(\theta)\triangleq {d\over
d\theta}\E\left[l(\bX(\theta))\cdot1_{\{h(\bX(\theta))\le
\xi\}}\right].
\end{equation}
One of the methods for estimating $\gamma'(\theta)$ is to explore
the possibility of interchanging the order of differentiation and
expectation, which is, unfortunately, invalid for the discontinuous
integrand in (\ref{eqn:DiscontTheta}). As a remedy, CMC suggests
choosing conditioning variables $\bY$, leading to
\begin{equation}\label{eqn:DiscontTheta1}
\gamma'(\theta)={d\over d\theta}\E\left[w(\bY;\theta)\right],
\end{equation}
where
\[w(\bY;\theta)\triangleq
\E\left[\left.l(\bX(\theta))\cdot1_{\{h(\bX(\theta))\le
\xi\}}\right|\bY\right].\]It is then expected that the interchange
of differentiation and expectation on (\ref{eqn:DiscontTheta1}) may
be valid for appropriate $\bY$, and if so, $\gamma'(\theta)$ can be
estimated by a sample mean of $dw(\bY;\theta)/d\theta$, provided
that it is easily computable.

A major limitation of CMC is on how to find appropriate conditioning
variables $\bY$. For some cases, finding $\bY$ to ensure both
differentiability of $w(\bY;\theta)$ and tractability of
$dw(\bY;\theta)/d\theta$ is difficult, e.g., when
$h(\bX)=\max(X_1,\dots,X_m)$.

In this paper, we study CMC from a change-of-variables perspective,
aiming to circumvent the difficulty on choosing $\bY$. Before we
proceed further, we decompose the problem by converting the
sensitivity w.r.t. any parameter $\theta$ to a sensitivity w.r.t.
$\xi$. In particular, Liu and Hong (2011, Theorem 1) showed that
under mild regularity conditions, $\gamma'(\theta)$ can be
represented as
\begin{equation}
\gamma'(\theta) = \E\left[\partial_\theta
l(\bX(\theta))\cdot1_{\{h(\bX(\theta))\le \xi\}}\right] -
\partial_{\xi} \E\left[l(\bX(\theta))\partial_\theta
h(\bX(\theta))\cdot1_{\{h(\bX(\theta))\le
\xi\}}\right],\label{eqn:LiuHong11}
\end{equation}
where $\partial$ denotes the differentiation operator, and
$\partial_\theta l(\bX(\theta))$ and $\partial_\theta
h(\bX(\theta))$ are pathwise derivatives that are often readily
computable from simulation. %For completeness, a detailed version of
%this result is provided in Section \ref{sec:LemmaLiuHong} of the
%appendix.

The result in (\ref{eqn:LiuHong11}) shows that the sensitivity
w.r.t. any $\theta$ can be related to the sensitivity w.r.t. $\xi$.
Note that the first term on the right-hand-side can be
straightforwardly estimated by a sample mean. The problem of
estimating $\gamma'(\theta)$ is then reduced to how to estimate the
second term.

Without loss of generality, the remainder of this paper is focused
on the reduced problem in which we want to estimate
\begin{equation}
\alpha'(\xi)\triangleq {d\over d\xi}\E\left[g(\bX)\cdot
1_{\{h(\bX)\le \xi\}}\right],\label{eqn:AlphaZ0}
\end{equation}
for some function $g$, where the dependence of $\bX$ on $\theta$ is
suppressed when there is no confusion.\footnote{There is no
smoothness requirement on the function $g$ which is allowed to be
discontinuous.}

\section{A Change-of-Variables Approach}\label{sec:Perspevtive}

This section introduces the main idea of the change-of-variables
approach. Section \ref{sec:Intuition} provides intuitions and
motivating examples. Section \ref{sec:Framework} lays down the
framework of the approach. Discussion of the approach in sensitivity
estimation for integrands of the form (\ref{eqn:IndForm}) will be
presented in a following section.

The change-of-variables argument applies to a general integrand,
beyond the form specified in (\ref{eqn:IndForm}). To avoid heavy
notations, in this section we work with a general integrand
$p(\bX)$, and discuss the change-of-variables approach for
$\E\left[p(\bX)\right]$.

\subsection{Intuitions and Motivating Examples}\label{sec:Intuition}

Consider the expectation $\E[p(\bX)]$. Recall that CMC yields
$\E[p(\bX)]=\E\left[w(\bY)\right]$, where
$w(\by)=\E\left[p(\bX)|\bY=\by\right]$. In terminology, we refer to
CMC based on conditioning variables as {\it conventional CMC}.

We highlight the intuition of our approach by recovering the above
conventional CMC estimator from a change-of-variables perspective,
where we ignore mathematical rigor at this stage. The
change-of-variables approach proceeds by constructing a one-to-one
mapping
\[u:\bX\mapsto (\bY,\bZ).\]Let $f(\bx)$ and $\tilde f(\by,\bz)$ denote the density functions of $\bX$
and $(\bY,\bZ)$, respectively. Then,
\begin{eqnarray}
\E[p(\bX)] = \int p(\bx)f(\bx)\,d\bx=\int\int
p(u^{-1}(\by,\bz))\tilde
f(\by,\bz)\,d\bz\,d\by.\label{eqn:MotiChange}
\end{eqnarray}
If we set
\[w(\by)=\left.\int p(u^{-1}(\by,\bz))\tilde f(\by,\bz)\,d\bz\right/\int \tilde
f(\by,\bz)\,d\bz,\]it can be easily verified that
\begin{eqnarray}
\E[p(\bX)] = \int w(\by)\int \tilde
f(\by,\bz)\,d\bz\,d\by=\E\left[w(\bY)\right],\label{eqn:MotiFubini}
\end{eqnarray}
where the second equality follows from the fact that $\int \tilde
f(\by,\bz)\,d\bz$ is the marginal density of $\bY$.

In short, the conventional CMC estimator can be recovered by a
change-of-variables argument in (\ref{eqn:MotiChange}) in
conjunction with an application of Fubini's Theorem in
(\ref{eqn:MotiFubini}). Compared to finding conditioning variables
$\bY$, the change-of-variables perspective offers more flexibility
in deriving CMC estimators. To see this, we consider two examples,
where $\bX=(X_1,X_2)$ with $X_1$ and $X_2$ following independent
exponential distributions with mean 1, and let $F_e$ and $f_e$
denote the c.d.f. and p.d.f. of the exponential distribution,
respectively.

\begin{exmp}\label{example:1}
\textnormal{$p(\bX)=X_1\cdot1_{\{X_1+X_2\le \xi\}}$ for a given
constant $\xi>0$. Conventional CMC may suggest choosing $X_1$ as a
conditioning variable. Then,
\[\E[p(\bX)] = \E\left(\E\left[X_1\cdot1_{\{X_1+X_2\le
\xi\}}|X_1\right]\right)=\E\left[X_1F_e(\xi-X_1)\right].\]This
estimator can be recovered from a change-of-variables perspective
using a one-to-one mapping: $(x_1,x_2)\mapsto
(x_1,x_1+x_2)\triangleq (y_1,y_2)$. The Jacobian (the absolute value
of the determinant of the derivative matrix) of the mapping is $1$,
and the density of $(Y_1,Y_2)\triangleq(X_1,X_1+X_2)$ is
$f_e(y_1)f_e(y_2-y_1)$. Then,
\[\E[p(\bX)] = \int\int 1_{\{y_2\le
\xi\}}y_1f_e(y_1)f_e(y_2-y_1)\,dy_2\,dy_1=\E\left[w(Y_1)\right],\]where
\[w(y_1) = {\int 1_{\{y_2\le \xi\}}y_1f_e(y_1)f_e(y_2-y_1)\,dy_2\over \int
f_e(y_1)f_e(y_2-y_1)\,dy_2}=y_1F_e(\xi-y_1).\]}

\textnormal{Interestingly, it is also possible to construct a
one-to-one mapping which leads to a new estimator that cannot be
derived by conventional CMC. Consider the one-to-one mapping:
\[(x_1,x_2)\mapsto (y_1,y_2,z)\triangleq \left({x_1\over x_1+x_2},{x_2\over
x_1+x_2},x_1+x_2\right).\]
\begin{figure}[H]\label{fig:Sum}
\centering
\includegraphics[width=0.4\textwidth]{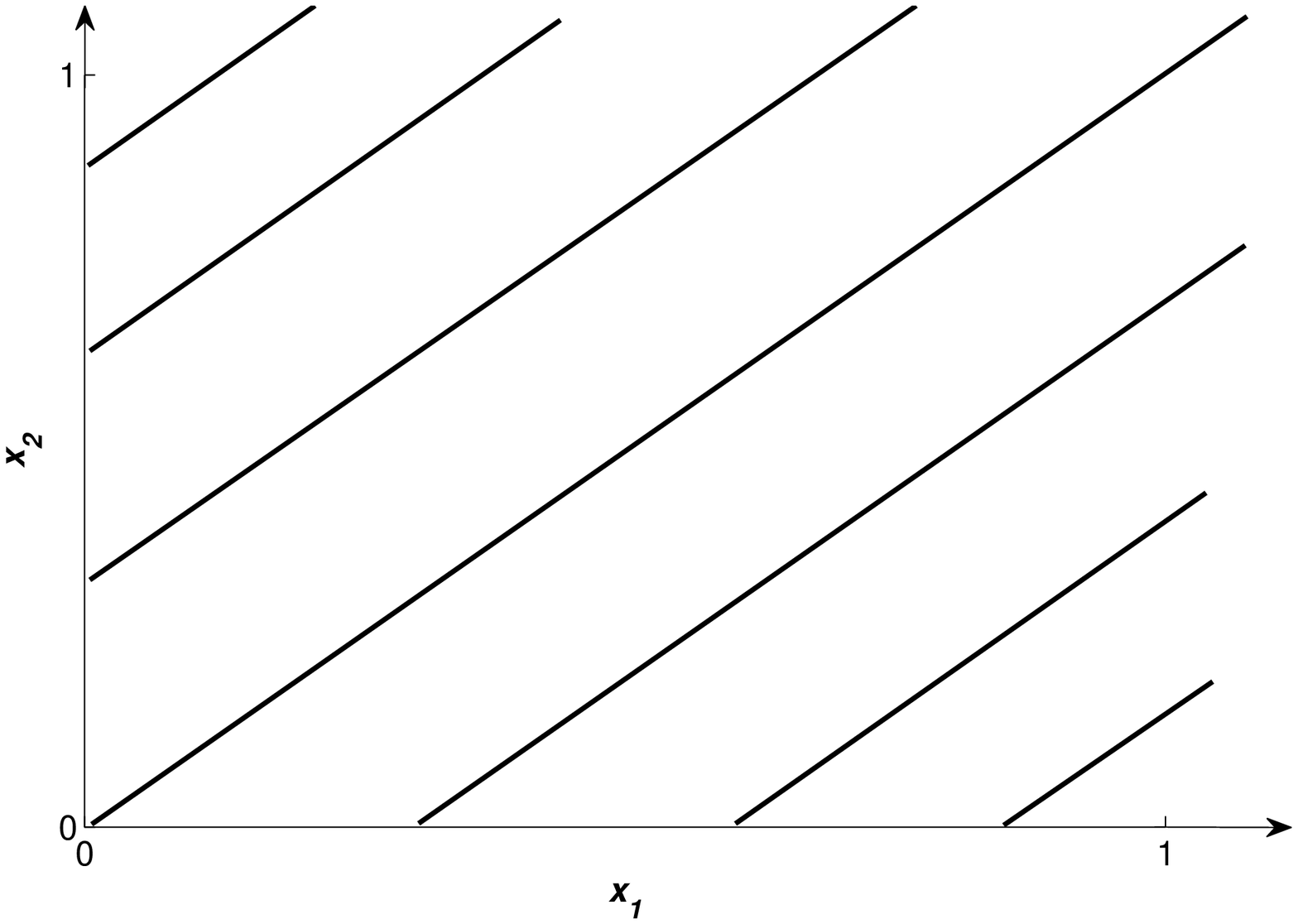}
\includegraphics[width=0.4\textwidth]{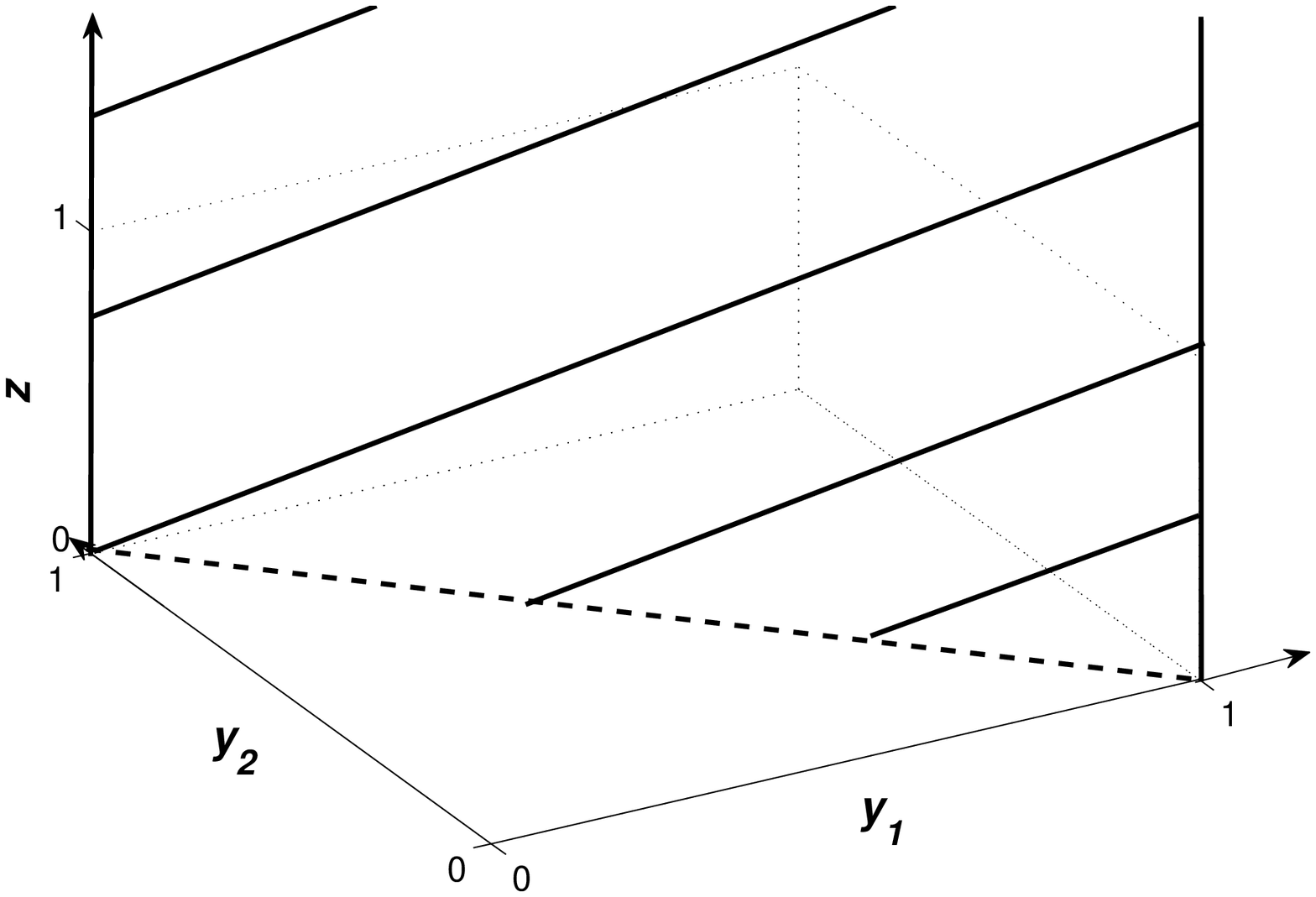}
\caption{Illustration of the mapping $(x_1,x_2)\mapsto
(y_1,y_2,z)\triangleq \left({x_1\over x_1+x_2},{x_2\over
x_1+x_2},x_1+x_2\right)$, where the left and right panels show the
domain and image sets of the mapping, respectively.}
\end{figure}
As shown in Figure \ref{fig:Sum}, image set of the mapping is
$\{(y_1,y_2,z): y_1+y_2=1,y_1>0,y_2>0,z>0\}$. It can be checked that
the mapping is indeed one-to-one, and the Jacobian of the mapping is
$\sqrt{2}/z$. Define $(Y_1,Y_2,Z)\triangleq \left({X_1\over
X_1+X_2},{X_2\over X_1+X_2},X_1+X_2\right)$, then the density of
$(Y_1,Y_2,Z)$ is $\tilde f(y_1,y_2,z)=zf_e(y_1z)f_e(y_2z)/\sqrt{2}$,
and
\begin{eqnarray}
\E[p(\bX)] = \int\int\int 1_{\{z\le \xi\}}y_1z\tilde
f(y_1,y_2,z)\,dz\,dy_1\,dy_2=\E\left[w(Y_1,Y_2)\right],\label{eqn:ExamTemp1}
\end{eqnarray}
where
\begin{align} \label{w_expression}% requires amsmath; align* for no eq. number
  w(y_1,y_2)={\int 1_{\{z\le \xi\}}y_1z\tilde f(y_1,y_2,z)\,dz\over \int\tilde
f(y_1,y_2,z)\,dz}={2y_1\over
y_1+y_2}-e^{-(y_1+y_2)\xi}\left({2y_1\over
y_1+y_2}-2y_1\xi-y_1(y_1+y_2)\xi\right).
\end{align}}
\end{exmp}

\begin{exmp}\label{example:2}
\textnormal{$p(\bX)=X_1\cdot 1_{\{\max(X_1,X_2)\le \xi\}}$.
Conventional CMC may suggest conditioning on $X_1$, leading to
\[\E[p(\bX)] = \E\left[F_e(\xi)X_11_{\{X_1\le \xi\}}\right].\]Similar to the argument in Example \ref{example:1}, this estimator follows from a change of variables with an identity
mapping.}

\textnormal{If we construct another one-to-one mapping:
\[(x_1,x_2)\mapsto (y_1,y_2,z)\triangleq \left({x_1\over \max(x_1,x_2)},{x_2\over
\max(x_1,x_2)},\max(x_1,x_2)\right).\]
%\begin{figure}[H]\label{fig:Max}
%\centering
%\includegraphics[width=0.4\textwidth]{broaden/original.eps}
%\includegraphics[width=0.4\textwidth]{broaden/max.eps}
%\caption{Illustration of the mapping $(x_1,x_2)\mapsto
%(y_1,y_2,z)\triangleq \left({x_1\over \max(x_1,x_2)},{x_2\over
%\max(x_1,x_2)},\max(x_1,x_2)\right)$, where the left and right
%panels represent the domain and image of the mapping respectively.}
%\end{figure}
Then by a similar argument as in Example \ref{example:1},
\begin{eqnarray*}
\E[p(\bX)] = \E\left[w(Y_1,Y_2)\right],
\end{eqnarray*}
where $(Y_1,Y_2)\triangleq \left({X_1\over \max(X_1,X_2)},{X_2\over
\max(X_1,X_2)}\right)$, and $w(y_1,y_2)$ has exactly the same form
as in (\ref{w_expression}).}
\end{exmp}

These examples illustrate that a change-of-variables argument offers
more flexibility for CMC and may lead to new estimators. However, it
should be pointed out that the integral in (\ref{eqn:ExamTemp1}) is
not well defined in Lebesgue sense. Note that the integral in
(\ref{eqn:ExamTemp1}) is taken over a set
$\{(y_1,y_2,z)|y_1>0,y_2>0,z>0,y_1+y_2=1\}$ which is a hyperplane in
$\mathbb{R}^3$. This set has Lebesgue measure 0 and the integral is
thus not well defined in Lebesgue sense. Rigorous definition of
integrals for such cases will be provided in Section
\ref{sec:Framework}.

%These two examples show that while conventional CMC estimators can
%be recovered by a change-of-variables argument, the alternative
%one-to-one mappings constructed above lead to new estimators that
%cannot be derived by conventional CMC. For Example 2, conventional
%CMC in the existing literature leads to estimators that are
%discontinuous in $\xi$. The alternative one-to-one mappings,
%however, yields an estimator $r(Y_1,Y_2)$ that is continuous in
%$\xi$. This continuity is particularly appealing in sensitivity
%estimation, which shall be made clear in Section
%\ref{sec:Sensitivity}.

%The derivation for the examples is based on constructing a
%one-to-one mapping and taking iterated integrations.  Mathematical
%treatment of this issue shall be provided in the following
%subsection.

\subsection{A Change-of-Variables Framework of
CMC}\label{sec:Framework}

We lay down a mathematical framework of the change-of-variables
approach for $\E\left[p(\bX)\right]$. Consider a one-to-one mapping
on $\Omega$, the support of $\bX$:
\[u\triangleq (u_1,\dots,u_n): \Omega\subset \mathbb{R}^m\mapsto \mathbb{R}^n,\]where $n\ge m$, $u_i$'s are functions of $\bx$, and the image set of
the mapping is denoted by $u\{\Omega\}\subset \mathbb{R}^n$. The
mapping is said to be differentiable at $\bx\in \mathbb{R}^{m}$ if
partial derivatives $Du(\bx)$ exist, and continuously differentiable
if all partial derivatives are continuous, where the partial
derivative matrix is defined as
\[ D u(\bx) = \left( \begin{array}{cccc}
\partial u_1(\bx)/\partial x_1 & \partial u_2(\bx)/\partial x_1 & \ldots & \partial u_n(\bx)/\partial x_1\\
\partial u_1(\bx)/\partial x_2 & \partial u_2(\bx)/\partial x_2 & \ldots & \partial u_n(\bx)/\partial x_2\\
\vdots & \vdots & \vdots & \vdots\\
\partial u_1(\bx)/\partial x_m & \partial u_2(\bx)/\partial x_m & \ldots & \partial u_n(\bx)/\partial x_m
\end{array} \right).\]Furthermore, for a subset $E\subset \mathbb{R}^{m}$, $u$ is said to be Lipschitz continuous if $|u(\bx_1)-u(\bx_2)|\leq
c|\bx_1-\bx_2|$ for any $\bx_1,\bx_2\in E$,
and a constant $c$. %The smallest constant $c$ such that the
%inequality holds for any $x,y\in E$ is denoted by
%$\mathrm{Lip}(u)=\mathrm{sup}\left\{|u(x)-u(y)|/|x-y|: x,y \in E,
%x\neq y\right\}$.

It should be emphasized that for a general one-to-one mapping, $n$
may not equal $m$. When $n>m$, the image set of the mapping is
essentially an ``$m$-dimensional" subset of $\mathbb{R}^n$. For
instance, the image set of the mapping illustrated in Figures
\ref{fig:Sum} is indeed a hyperplane that is a 2-dimensional subsets
of $\mathbb{R}^3$. For such cases, the integrals taken over the
image sets are not well defined in Lebesgue sense. To deal with this
issue, we work with Hausdorff measure, an extension of Lebesgue
measure.

\begin{definition}[Hausdorff Measure on $\mathbb{R}^n$]
Let $A$ be a nonempty subset of $\mathbb{R}^n$, and define its
diameter by $\text{diam}(A) \triangleq\{|x-y|: x,y\in A\}$. For any
$\delta\in(0,\infty]$ and $s\in [0,\infty)$, define $\mathcal
{H}^s_\delta(A)= \inf
\left\{\sum_{i=1}^{\infty}\omega(s)\left({\text{diam}(A_i)\over
2}\right)^s; A\subset \cup_{i=1}^\infty A_i, \text{diam}(A_i)\le
\delta\right\}$, where $\omega(s)= \pi^{s/2}/\Gamma(s/2+1)$, and
$\Gamma(t)=\int_{0}^{\infty}e^{-x}x^{t-1}\,dx$. Then $s$-dimensional
Hausdorff measure on $\mathbb{R}^n$ is defined as
\[\mathcal {H}^s(A)=\lim_{\delta\rightarrow 0}\mathcal
{H}^s_\delta(A).\]Furthermore, Hausdorff dimension of $A$ is defined
by $\dim_{\mathcal H}(A)=\inf \left\{s\in [0,\infty);{\mathcal
H}^s(A)=0\right\}$.
%\[\dim_{\mathcal H}(A)=\inf \left\{s\in [0,\infty);{\mathcal H}^s(A)=0\right\}.\]
\end{definition}

%\begin{remark}
%The dimension $\dim_{\mathcal H}(A)$ need not be an integer. Even if
%$\dim_{\mathcal H}(A)=k$ os an integer and $0<{\cal H}^k(A)<\infty$,
%$A$ need not be a ``$k$-dimensional surface" in any sense; see
%examples of extremely complicated Cantor-like subsets.
%\end{remark}

Essentially, Hausdorff measure generalizes the concepts of length,
area, and volume. For instance, $1$-dimensional Hausdorff measure of
a smooth curve in $\mathbb{R}^n$ is the length of the curve, and
$2$-dimensional Hausdorff measure of a smooth surface in
$\mathbb{R}^n$ is its area. Hausdorff measure has several
properties: ${\mathcal H}^s \equiv 0$ on $\mathbb{R}^n$ for any
$s>n$, ${\mathcal H}^s(\lambda A)=\lambda^s{\mathcal H}^s(A)$ for
any $\lambda>0$, and ${\mathcal H}^s$ is equivalent to Lebesgue
measure on $\mathbb{R}^n$ when $s=n$. Hausdorff measure is a useful
tool to measure ``small" subsets in $\mathbb{R}^n$. For instance, a
set $A=[0,1]\times [0,1]\times \{1\}\subset \mathbb{R}^3$ is a
``copy" of the 2-dimensional square $[0,1]\times [0,1]$ in
$\mathbb{R}^3$. Its volume is 0, as its Lebesgue measure is $0$. It
has, however, an area of 1, i.e., ${\cal H}^2(A)=1$, and its
Hausdorff dimension is 2.

%For any one-to-one mapping $u:\mathbb{R}\mapsto \mathbb{R}^n$ which
%is differentiable. For given constants $a$ and $b$, the set
%$u\{[a,b]\}\subset \mathbb{R}^n$ defines a curve in $\mathbb{R}^n$.
%The Hausdorff measure of this set is then the length of the curve,
%which can be calculated by $\int_a^b |\partial u(t)|\,dt$, where
%$\partial$ dentes the differentiation operator.

%\begin{definition}[Continuously Differentiable ($\mathcal {C}^1$) Mappings]\label{def:C1}
%A mapping $u:\mathbb{R}^{m}\mapsto \mathbb{R}^{n}$ is differentiable
%at $x\in \mathbb{R}^{m}$ if the partial derivatives $Du(x)$ exist.
%Furthermore, $u$ is continuously differentiable if all partial
%derivatives are continuous.
%\end{definition}

To characterize smoothness of a set on $\mathbb{R}^n$, we introduce
the concept of rectifiability.

\begin{definition}[Rectifiable Sets]\label{def:rectifiable}
For $m\le n$, a set $A\subset\mathbb{R}^n$ is said to be
$(\mathcal{H}^{m},m)$ rectifiable if there exist a countable
collection $\{v_i\}_{i\ge 1}$ of continuously differentiable
mappings $v_i:\mathbb{R}^m \mapsto \mathbb{R}^n$ such that ${\cal
H}^m\left(A\backslash \bigcup_{i=1}^\infty
v_i\left(\mathbb{R}^m\right)\right)=0$.
\end{definition}

Typically, a set is rectifiable if it can be represented as the
union of images of continuously differentiable mappings. For
instance, the set $A=[0,1]\times [0,1]\times \{1\}\subset
\mathbb{R}^3$ is $({ \cal H}^2,2)$ rectifiable, because it is the
image of a mapping $u:\mathbb{R}^2\mapsto \mathbb{R}^3$ where
$u(x_1,x_2)=(x_1,x_2,1)$. Intuitively speaking, a rectifiable set is
a piece-wise ``smooth" set in $\mathbb{R}^n$. It has many of the
desirable properties when taking Hausdorff integrals.

When applying changes of variables for the mapping $u:
\mathbb{R}^m\mapsto \mathbb{R}^n$, the Jacobian plays a critical
role. When $n=m$, the Jacobian is defined as the absolute value of
the determinant of $Du$, serving as a corrective factor that relates
the ``volumes" of the domain and image sets. In a more general
setting where $n\ge m$, $Du$ may not be a square matrix, and an
extension of the Jacobian is defined as follows(see, e.g., Morgan
(2009)).

\begin{definition}[$k$-dimensional Jacobian]
Let $J_{k}^u(\bx)$ denote the $k$-dimensional Jacobian of $u$ at
$\bx$. Then $J_k^u(\bx)=0$ if $\text{rank}(Du(\bx))<k$, and it is
equal to square root of the sum of squares of the determinants of
the $k\times k$ submatrices of $Du(\bx)$ if $\text{rank}(Du(\bx))\ge
k$.
\end{definition}
%\begin{definition}[$k$-dimensional Jacobian]
%Geometrically, $J_{k}(u(a))$ is the maximum $k$-dimensional volume
%of the image under $Du(a)$ of a unit $m$-dimensional cube in
%$\mathbb{R}^{n}$. In particular, if rank $Du(a)<k$, $J_{k}(u(a))=0$.
%If rank $Du(a)<k$, then the square of Jacobian $J_{k}^{2}(u(a))$
%equals the sum of the squares of the determinants of the $k\times k$
%submatrices of $Du(a)$.
%\end{definition}
%The definition implies that when the rank of the derivative matrix
%equals $m$, then
%$J_{m}^{2}(u(a))=\sqrt{\mathrm{det}(Du(a)(Du(a))^{T})}$. In a
%special case when $m=n$, then $J_{m}(u(a))$ is the absolute value of
%the determinant of $Du(a)$, which is the usual Jacobian used in
%elementary calculus.

%\begin{definition}[Lipschitz Continuous Mappings]\label{Def:lips}
%Let $E\subset \mathbb{R}^{m}$ be a subset. A mapping $u:E\mapsto
%\mathbb{R}^{n}$ is Lipschitz if $|u(x)-u(y)|\leq C|x-y|$ for any
%$x,y\in E$, for a constant C. The smallest constant $C$ such that
%the inequality holds for any $x,y\in E$ is denoted by
%$\mathrm{Lip}(u)=\mathrm{sup}\left\{\frac{|u(x)-u(y)|}{|x-y|}; x,y
%\in E, x\neq y\right\}$.
%\end{definition}

In what follows, we discuss how the change-of-variables approach
works as in Examples \ref{example:1} and \ref{example:2}. It
proceeds by two steps. In the first step, we note that $m\le n$,
$u\{\Omega\}$ is essentially an $m$-dimensional subset on
$\mathbb{R}^n$ because $u$ is one-to-one. Then
$\E\left[p(\bX)\right]$, a Lebesgue integral taken over $\Omega$, is
equal to a Hausdorff integral taken over the set $u\{\Omega\}$,
which is justified by a change-of-variables formula for Hausdorff
measure.

The second step is to take iterated integration by using an extended
version of Fubini's Theorem. To this end, we first note that
$u\{\Omega\}$ can be written in a product form of
$u\{\Omega\}=S_1\times S_2$ under mild smoothness conditions, where
$S_1\subset \mathbb{R}^{n-q}$ and $S_2\subset\mathbb{R}^q$ for some
integer $q\le n$. Then similar to that in Example \ref{example:1},
taking integral\footnote{ Hausdorff dimension of $S_1$ plays an
important role in taking integrals. If Hausdorff dimension of $S_1$
is $k$ for some $k\le q$, Hausdorff dimension of $S_2$ should be
$m-k$, because Hausdorff dimension of $u\{\Omega\}$ is $m$.} over
the set $S_1$ yields a function of $w$, which is in fact a function
of the elements in $S_2\subset \mathbb{R}^q$. More generally, this
step can be viewed through an argument on mappings. Specifically,
construct a mapping
\[\varphi: u\{\Omega\}\mapsto \mathbb{R}^q.\]
Suppose Hausdorff dimension of the image set
$\varphi\{u\{\Omega\}\}$ is $k$ for some $k\le q$. Then for any
fixed $\bz\in \varphi\{u\{\Omega\}\}$, Hausdorff dimension of the
set $\varphi^{-1}\{\bz\}$ is $m-k$. Then fixing $\bz$ and taking
integral over $\varphi^{-1}\{\bz\}$ lead to a function of $w$ as in
Examples 1 and 2. This result is stated in the following theorem,
whose proof is provided in the appendix.

\begin{theorem}\label{thm: expecttransform}
Assume that Assumption \ref{assum:Density} holds. Suppose $n\ge m\ge
k\le q$, $u:\Omega\mapsto \mathbb{R}^n$ is a continuously
differentiable one-to-one mapping, and $\varphi: u\{\Omega\}\mapsto
\mathbb{R}^q$ is Lipschitz continuous. If $\varphi\{u\{\Omega\}\}$
is $({\cal H}^{k},k)$ rectifiable, and $J_{m}^u(\bx)\neq 0$ and
$J_{k}^\varphi(\by)\neq 0$ for almost all $\bx\in \Omega$ and
$\by\in u\{\Omega\}$ up to null sets, then
\begin{equation}\label{transformation}
\E\left[p(\bX)\right]=\E\left[w\left(\varphi(u(\bX))\right)\right],
\end{equation}
where the function $w$ is defined by
\begin{equation}
w(\bz)\triangleq
\left.\int_{\varphi^{-1}\{\bz\}}\frac{p(u^{-1}(\by))f(u^{-1}(\by))}{J_{m}^u(u^{-1}(\by))J_{k}^\varphi(\by)}\mathcal{H}^{m-k}(d\by)\right/\int_{\varphi^{-1}\{\bz\}}\frac{f(u^{-1}(\by))}{J_{m}^u(u^{-1}(\by))J_{k}^\varphi(\by)}\mathcal{H}^{m-k}(d\by).\label{eqn:Weight}
\end{equation}
Moreover, $w\left(\varphi(u(\bX))\right)$ has a smaller variance
than $p(\bX)$, i.e.,
\[\Var\left[w\left(\varphi(u(\bX))\right)\right]\le \Var\left[p(\bX)\right].\]
\end{theorem}

Theorem \ref{thm: expecttransform} shows that $w(\varphi(u(\bX)))$
has the same expectation as $p(\bX)$, and provides a closed-form
expression for the function $w$. In addition, it shows that
$w(\varphi(u(\bX)))$ has a smaller variance than $p(\bX)$, which is
expected as in conventional CMC.

Theorem \ref{thm: expecttransform} works with Hausdorff integrals,
aiming to accommodate a more general setting. For a special case
when $n=m$, $k=q<n$, and the mapping $\varphi$ is specified by
$\varphi(\by)=(y_1,\dots,y_k)$ for $\by\triangleq (y_1,\dots,y_m)$,
Hausdorff integrals are equivalent to Lebesgue integrals, and
$w\left(\varphi(u(\bX))\right)$ is equivalent to
$\E\left[p(\bX)|u_1(\bX),\dots,u_k(\bX)\right]$.

%The second $w$ function in Example 1 can be derived by applying
%Theorem 1 with $n=m+1$, $k=m-1$, $q=m$, and
%\[u(\bx)=(\bx/h(\bx),h(\bx)),\quad \varphi(\by)=(y_1,\dots,y_m), \forall \ \bx\triangleq
%(x_1,\dots,x_m),\ \by\triangleq (y_1,\dots,y_{m+1}).\]

To better understand the insights offered by Theorem \ref{thm:
expecttransform}, we look into the result in greater detail. Note
that by the theorem,
\[\E\left[w\left(\varphi(u(\bX))\right)\right]=\E\left[p(\bX)\right]=\E\left(\E\left[p(\bX)|\varphi(u(\bX))\right]\right).\]It
is then reasonable to expect that
$w\left(\varphi(u(\bX))\right)=\E\left[p(\bX)|\varphi(u(\bX))\right]$,
which is also confirmed as a by-product in the proof of the theorem.
Compared to conventional CMC that directly conditions on
$\varphi(u(\bX))$, Theorem \ref{thm: expecttransform} is meaningful
in the following two aspects.

First, Theorem \ref{thm: expecttransform} provides a closed-form
expression of $\E\left[p(\bX)|\varphi(u(\bX))\right]$, that applies
to a very general simulation model of $\bX$, in contrast to
conventional CMC where deriving explicit formula of the conditional
expectation often relies highly on the structure of $\bX$. As seen
in Examples 1 and 2, the expression of the function $w$ may lead to
estimators with explicit formulas. Derivation of such explicit
formulas shall be further exemplified in Section
\ref{sec:Application} for applications under more practical
settings.

Second, Theorem \ref{thm: expecttransform} offers great flexibility
for us to choose among various mappings $u$ and $\varphi$, so as to
make use of the form of the integrand $p(\bX)$. This point shall be
made clearer in Section \ref{sec:Sensitivity} where we construct
mappings $u$ and $\varphi$ to exploit the unique structure of a
discontinuous integrand of the form (\ref{eqn:IndForm}), leading to
efficient sensitivity estimators.

%the availability of the expression $w(\cdot)$ may ease the
%estimation of $\E\left[p(\bX)|\varphi(u(\bX))\right]$ for cases
%where no explicit formulas can be derived. For such cases, one may
%resort to numerical or Monte Carlo integration methods to
%approximate or estimate the conditional expectation for any given
%realization of $\varphi(u(\bX))$. However, numerical or Monte Carlo
%integration could be difficult for some cases. For instance, when
%$\varphi(u(\bX))=\bX/\sum_{i=1}^{m}X_i$ as in Example 1, it is not
%clear how to estimate the conditional expectation directly without
%using $w(\cdot)$, while both the numerator and denominator of $w$
%are indeed one-dimensional Lebesgue integrations (see Section
%\ref{sec:Sensitivity}) and can be easily approximated by numerical
%integration methods.

%Last but not least, the expression $w(\cdot)$ is particularly
%appealing when the integrand $p(\bX)$ is discontinuous and one is
%interested in estimating the sensitivity of $\alpha$. For such
%cases, a main goal CMC is to integrate out the discontinuities of
%the integrand and ensure the conditional expectation is continuous.
%To this end, $w(\cdot)$ is appealing in that the mappings $\varphi$
%and $u$ can be chosen appropriately such that the integrations in
%(\ref{eqn:Weight}) are taken along the space where discontinuities
%occur, thus leading to an expression of $w(\varphi(u(\bX)))$ that is
%naturally continuous. A deeper investigation of this issue is the
%focus of the following section.

\section{Change-of-Variables Approach for Sensitivity
Estimation}\label{sec:Sensitivity}

In this section we study the change-of-variables approach for
sensitivity estimation for a discontinuous integrand of the form
(\ref{eqn:IndForm}). In particular, we first provide a general
estimator in Theorem \ref{thm:Sensitivity} for the sensitivity,
which applies to general mappings that satisfy mild conditions. We
then consider two cases in Sections \ref{sec:Case1} and
\ref{sec:Case2}, providing sensitivity estimators when the function
$h$ in (\ref{eqn:IndForm}) satisfies Conditions \ref{con:C1} and
\ref{con:C2}, respectively. These estimators require mainly
smoothness conditions on $h$, and their implementation do not
require any construction of mappings. We find that the function $h$
encountered in many practical applications satisfies Condition
\ref{con:C1} and/or Condition \ref{con:C2}. The proposed estimators
are, therefore, useful for a wide range of applications.

Recall that our problem is to estimate
\begin{equation}\label{eqn:DefSensitivity}
\alpha'(\xi)\triangleq {d\over d\xi}\E\left[g(\bX)\cdot
1_{\{h(\bX)\le \xi\}}\right],
\end{equation}
for given functions $g$ and $h$, as stated in Section
\ref{sec:Problem}.

To facilitate analysis, we make the following smoothness assumption
on the function $h$. Moreover, we let $(c_0,c_1)$ denote the
interval in which $h(\bx)$ may take values for $\bx\in\Omega$, where
$c_0=-\infty$ and $c_1=\infty$ are allowed. To avoid trivial cases,
we assume $\xi\in (c_0,c_1)$.
\begin{assumption}\label{assum:h}
The function $h(\bx)$ is continuously differentiable for almost all
$\bx\in \Omega$ up to a Lebesgue null set.
\end{assumption}

To analyze the estimation of $\alpha'(\xi)$, we start our discussion
with the result in Theorem \ref{thm: expecttransform}, which can be
applied to obtain a function of $w$ such that
\[\alpha(\xi)=\E\left[g(\bX)\cdot
1_{\{h(\bX)\le
\xi\}}\right]=\E\left[w\left(\varphi(u(\bX));\xi\right)\right],\]where
$\xi$ is added as an argument of the function $w$ to explicitly
account for the dependence. To enable the use of IPA, one may
attempt to construct mappings $u$ and $\varphi$ such that
$w\left(\varphi(u(\bX));\xi\right)$ is smooth in $\xi$. To this end,
it is natural to take into account the unique structure of
$g(\bX)\cdot 1_{\{h(\bX)\le \xi\}}$ that $\xi$ appears in the
indicator function only.

Inspired by Examples \ref{example:1} and \ref{example:2}, we
consider the following mappings:
\begin{equation}\label{eqn:mappings}
u: \bx\in \Omega\mapsto (u_1(\bx),\dots,u_m(\bx),h(\bx)), \quad
\varphi: \by\in u\{\Omega\}\mapsto (y_1,\dots,y_m),
\end{equation}
for $\bx\triangleq (x_1,\dots,x_m)$ and $\by\triangleq
(y_1,\dots,y_{m+1})$, where $\{u_i, 1\le i\le m\}$ are continuously
differentiable functions of $\bx$.

A common feature of the $u$ mappings we consider is that $h(\bx)$
serves as the last dimension in the image set, taking integral along
which removes discontinuity in the integrand $g(\bx)\cdot
1_{\{h(\bx)\le \xi\}}$. The resulting
$w\left(\varphi(u(\bX));\xi\right)$ is thus continuous in $\xi$,
which enables the use of IPA. For the sake of using
change-of-variables formula, we require that $u$ is one-to-one.

It should also be pointed out that $u$ in (\ref{eqn:mappings}) is in
general a mapping from $\Omega\subset \mathbb{R}^m$ to
$\mathbb{R}^{m+1}$. To ensure that it is one-to-one, Hausforff
dimension of the image set $u\{\Omega\}$ should be $m$. In other
words, $u\{\Omega\}$ is an $m$-dimensional subset in
$\mathbb{R}^{m+1}$. Furthermore, for the integral taken along the
last dimension to make sense, Hausdorff dimension of the last
coordinate of $u\{\Omega\}$ should be $1$, which intuitively implies
that for fixed values of the first $m$ coordinates, the last
dimension of $u\{\Omega\}$ contains at least a closed interval. In
plain words, for any fixed $\bz\in \mathbb{R}^m$, the curve
$\{t\in\mathbb{R}: t=h(\bx), (u_1(\bx),\dots,u_m(\bx))=\bz\}$ has a
positive length. Mathematically, this requirement can be translated
to
\[{\cal H}^1\left(\varphi^{-1}\{\bz\}\right)>0,\quad \forall\ \bz\in
\varphi\{u\{\Omega\}\},\]where by the definition of $\varphi$, the
set $\varphi^{-1}\{\bz\}$ is
\[\varphi^{-1}\{\bz\}=\{(\bz,t)\in u\{\Omega\}:t\in
(c_0,c_1)\}.\]

%To fulfill the conditions of Theorem \ref{thm: expecttransform}, we
%require $u$ to be a continuously differentiable one-to-one mapping.
%In other words, the functions $h(\cdot)$ and $u_i(\cdot)$,
%$i=1,\dots,m$ are continuously differentiable, and $u\{\Omega\}$ is
%a subset of $\mathbb{R}^{m+1}$ that has Hausdorff dimension $m$.

%The reason why the mappings in (\ref{eqn:mappings}) are appealing
%can be seen intuitively as follows. Recall that a main goal of
%constructing the mapping is to ensure the $w(\cdot)$ as in
%(\ref{eqn:Weight}) is continuous in $\xi$ that appears only in the
%numerator of $w(\cdot)$ through the function $p(\bx)=g(\bx)\cdot
%1_{\{h(\bx)\le \xi\}}$. Take a closer look at the integration in the
%numerator of $w(\cdot)$, which is taken over a set
%$\varphi^{-1}\{\bz\}$ for any fixed $\bz\in
%\varphi\{u\{\Omega\}\}\subset\mathbb{R}^m$. Fixing $\bz$, the set
%$\varphi^{-1}\{z\}$ is in fact a subset spanned by the last
%coordinate while fixing the first $m$ coordinates. Therefore, the
%indicator $1_{\{h\le \xi\}}$ inside the integration can be replaced
%by integration limits along $h$ over $(-\infty,\xi)$, and then it
%can be seen that the integration is continuous in $\xi$.

Given these conditions, we apply Theorem \ref{thm: expecttransform}
and have
\begin{eqnarray*}
\alpha(\xi)=\E\left[g(\bX)\cdot 1_{\{h(\bX)\le
\xi\}}\right]=\E\left[w\left(\varphi(u(\bX));\xi\right)\right],
\end{eqnarray*}
where
\begin{eqnarray}
w(\bz;\xi) =
\left.\int_{c_0}^{\xi}{g(u^{-1}(\bz,t))f(u^{-1}(\bz,t))\over
J_m^u(u^{-1}(\bz,t))}dt\right/\int_{c_0}^{c_1}{f(u^{-1}(\bz,t))\over
J_m^u(u^{-1}(\bz,t))}dt,\label{eqn:SenProofTemp1}
\end{eqnarray}
where the Jacobians of $\varphi$ have been cancelled out because
$J_{m-1}^\varphi(\by)\equiv \sqrt{m}$ for all $\by$.

The function $w(\bz;\xi)$ in (\ref{eqn:SenProofTemp1}) is typically
continuous in $\xi$, because $\xi$ is simply a limit of integration.
It is then reasonable to expect that
\[\alpha'(\xi)={d\over d\xi}\E\left[w\left(\varphi(u(\bX));\xi\right)\right]=\E\left[{d\over
d\xi}w\left(\varphi(u(\bX));\xi\right)\right].\]This result is
summarized in the following theorem, whose proof is provided in the
appendix.

\begin{theorem}\label{thm:Sensitivity}
Suppose that the mappings $u$ and $\varphi$ are specified as in
(\ref{eqn:mappings}), and $u$ is one-to-one and continuously
differentiable. If Assumptions \ref{assum:Density}-\ref{assum:h}
hold, $J_m^u(\bx)\neq 0$ for almost all $\bx\in \Omega$, and ${\cal
H}^1\left(\varphi^{-1}\{\bz\}\right)>0$ for almost all $\bz\in
\varphi\{u\{\Omega\}\}$, then
\[\alpha(\xi)=
\E\left[w\left(\varphi(u(\bX));\xi\right)\right],\]where
$w(\cdot;\xi)$ is defined in (\ref{eqn:SenProofTemp1}).

If, in addition, $w\left(\varphi(u(\bX));\cdot\right)$ satisfies a
Lipschitz continuity condition w.r.t. $\xi$, i.e., there exists a
random variable $K$ with $\E(K)<\infty$ such that for all small
enough $\Delta$,
\[|w\left(\varphi(u(\bX));\xi+\Delta\right)-w\left(\varphi(u(\bX));\xi\right)|\le
K|\Delta|,\quad \text{with probability 1},\]then,
\[\alpha'(\xi) =
\E\left[\nu(\varphi(u(\bX));\xi)\right],\]where
\[\nu(\bz;\xi)=\left.{g(u^{-1}(\bz,\xi))f(u^{-1}(\bz,\xi))\over
J_m^u(u^{-1}(\bz,\xi))}\right/\int_{c_0}^{c_1}{f(u^{-1}(\bz,t))\over
J_m^u(u^{-1}(\bz,t))}dt.\]
\end{theorem}

\begin{remark}
\textnormal{A special case of the result of Theorem
\ref{thm:Sensitivity} is $g(\bx)\equiv 1$ in which the quantity of
interest, $\alpha'(\xi)$, is in fact the density of $h(\bX)$
evaluated at $\xi$. When Theorem \ref{thm:Sensitivity} is
applicable, a sample-mean estimator can be derived for the density
of $h(\bX)$.}
\end{remark}

Theorem \ref{thm:Sensitivity} shows that the sensitivity
$\alpha'(\xi)$ is equal to $\E\left[
\nu(\varphi(u(\bX));\xi)\right]$, and provides a closed-form
expression of the function $\nu$. When the one-dimensional integral
in $\nu$ can be derived in explicit forms, an explicit formula of
$\nu$ can be obtained, as in Examples \ref{example:1} and
\ref{example:2}. When explicit formulas are not available, one may
resort to one-dimensional numerical integration methods to
approximate $\nu(\varphi(u(\bX));\xi)$ for any realization of
$\varphi(u(\bX))$. Efficient one-dimensional numerical integration
tools are available in many commercial softwares such as Matlab.

To ensure the validity of interchanging differentiation and
expectation, Theorem \ref{thm:Sensitivity} requires that
$w\left(\varphi(u(\bX));\xi\right)$ satisfies a Lipschitz continuity
condition w.r.t. $\xi$, which has been a commonly used condition in
sensitivity estimation literature; see, e.g., Broadie and Glasserman
(1996) and Liu and Hong (2011). This assumption typically does not
impose any obstacle in practice, because the interchange is usually
valid when $w\left(\varphi(u(\bX));\xi\right)$ is continuous in
$\xi$, which is obviously true by the way $w$ is defined.

A key condition of Theorem \ref{thm:Sensitivity} is ${\cal
H}^1\left(\varphi^{-1}\{\bz\}\right)>0$ for any given $\bz$, which
may not hold in general. This condition also depends on how the $u$
and $\varphi$ mappings are constructed. In practice, before applying
the result in Theorem \ref{thm:Sensitivity}, one may need to
construct the mappings and verify this condition accordingly. To
reduce simulation practitioners' effort in carrying out
verifications, Sections \ref{sec:Case1} and \ref{sec:Case2} consider
two cases for which we provide two sufficient conditions imposed on
the function $h$ only, and specify explicitly what mappings to be
used. These two sufficient conditions are easier to verify. We also
note that the function $h$ in many practical applications satisfies
either one of the sufficient conditions, or both.

\subsection{Case 1}\label{sec:Case1}

Consider a function $h$ that satisfies the following condition:
\begin{condition}\label{con:C1}
The mapping $\bx\mapsto (x_1,\dots,x_{m-1},h(\bx))$ is one-to-one,
and $\Pr\{\partial_m h(\bX)=0\}=0$.
\end{condition}
We claim that when Condition \ref{con:C1} holds, the key condition
of Theorem \ref{thm:Sensitivity}, i.e., ${\cal
H}^1\left(\varphi^{-1}\{\bz\}\right)>0$ for any given $\bz$, is
satisfied for appropriate mappings $u$ and $\varphi$. To see this,
we set $u_i(\bx)=x_i$ for $1\le i\le m-1$ and $u_m(\bx)\equiv 1$ in
(\ref{eqn:mappings}). Note that $u_m(\bx)$ takes a constant value.
The mappings in (\ref{eqn:mappings}) are therefore equivalent to the
following:
\begin{equation}\label{eqn:MappingSimple}
\tilde u: \bx \mapsto(x_1,\dots,x_{m-1},h(\bx)),\quad \tilde
\varphi:\by\mapsto (y_1,\dots,y_{m-1}),
\end{equation}
where $\by\triangleq (y_1,\dots,y_{m})$.

To see why Condition \ref{con:C1} implies ${\cal
H}^1\left(\tilde\varphi^{-1}\{\bz\}\right)>0$ for any given $\bz$,
we note that when Condition \ref{con:C1} is satisfied, the set
$\tilde \varphi^{-1}\{\bz\}$ for a given $\bz$ is
$\{(\bz,t):t\in\mathbb{R}, t=h(\bx)\ \text{for some}\ \bx\in
\Omega\}$, which is indeed a straight line with a positive length
because $\Omega$ has at least a compact subset and $h$ is continuous
by Assumption \ref{assum:h}. The condition that ${\cal
H}^1\left(\tilde \varphi^{-1}\{\bz\}\right)>0$ for almost all $\bz$
is therefore satisfied.

It can be easily seen that the Jacobian of $\tilde u$ is
\[\tilde J(\bx)=|\partial_m h(\bx)|,\]where $\partial_m h(\bx)$ denotes
the partial derivative of $h(\bx)$ w.r.t. $x_m$. Denote
$(y_1,\dots,y_{m-1})$ by $\bz$. For any given $\by=(\bz,y_m)$, let
$v(\bz,y_m)$ denote the solution to the equation $h(\bz,x_m)=y_m$,
which is unique because $\tilde u$ is one-to-one. Then it can be
easily seen that
\[\tilde u^{-1}(\by)=\left(\bz,v(\bz,y_m)\right),\]and $h\left(\tilde
u^{-1}(\by)\right)=y_m$.

Applying Theorem \ref{thm:Sensitivity}, we obtain an expression of
$\alpha'(\xi)$. This result is summarized in the following
proposition, whose proof is straightforward and thus omitted.

\begin{pro}\label{pro:Simple}
Suppose Assumptions \ref{assum:Density}-\ref{assum:h} and Condition
\ref{con:C1} are satisfied. Then,
\[\alpha(\xi)=\E\left[\widetilde
w\left(\bX_{-m};\xi\right)\right],\]where $\bX_{-m}\triangleq
(X_1,\dots,X_{m-1})$, and
\[\widetilde w(\bz;\xi)=\left.\int_{c_0}^{\xi}{g\left(\bz,v(\bz,t)\right)f(\bz,v(\bz,t))\over \tilde J(\bz,v(\bz,t))}\,dt\right/\int_{c_0}^{c_1}{f(\bz,v(\bz,t))\over \tilde J(\bz,v(\bz,t))}\,dt.\]
If, in addition, $\widetilde w(\bX_{-m};\cdot)$ satisfies a
Lipschitz continuity condition w.r.t. $\xi$ as in Theorem
\ref{thm:Sensitivity}, then,
\[\alpha'(\xi)=\E\left[\widetilde
\nu\left(\bX_{-m};\xi\right)\right],\]where
\begin{equation}\label{eqn:EstCase1}
\widetilde
\nu(\bz;\xi)=\left.{g\left(\bz,v(\bz,\xi)\right)f(\bz,v(\bz,\xi))\over
\tilde
J(\bz,v(\bz,\xi))}\right/\int_{c_0}^{c_1}{f(\bz,v(\bz,t))\over
\tilde J(\bz,v(\bz,t))}\,dt.
\end{equation}
\end{pro}

\begin{remark}
\textnormal{It should be pointed out that the result in Proposition
\ref{pro:Simple} still applies if we change the mapping in Condition
\ref{con:C1} to
\[\bx\mapsto \left(x_{i_1},\dots,x_{i_{m-1}},h(\bx)\right)\]for any
permutation $(i_1,\dots,i_m)$ of $\{1,\dots,m\}$. We only need to
change the indices of the function arguments accordingly. This
offers more flexibility to verify Condition \ref{con:C1} in
practice.}
\end{remark}

Proposition \ref{pro:Simple} shows that under appropriate
conditions, $\alpha'(\xi)$ can be estimated by a sample mean of
$\widetilde \nu\left(\bX_{-m};\xi\right)$, where the function
$\widetilde \nu$ is specified in (\ref{eqn:EstCase1}). For any given
$\bz$, evaluation of the function $\widetilde \nu$ requires knowing
$v(\bz,t)$ that is in fact the inverse of $h(\bz,x_m)$ as a function
of $x_m$. When an explicit formula of $v(\bz,t)$ can be derived, the
evaluation is straightforward. When explicit formulas are not
available, in principle it can be approximated by numerically
solving the equation $h(\bz,x_m)=t$ by using, e.g., one-dimensional
line search methods. For instance, when $t=\xi$, $v(\bz,\xi)$ can be
approximated by line search algorithms for any given $\bz$. However,
evaluation of the denominator term of $\widetilde \nu$ may require
line search operations for all possible $\bz$'s and a grid of $t$'s
in $(c_0,c_1)$, which may be too computationally expensive to
afford.

To resolve this issue, we may further apply a change of variables on
the denominator term. Specifically, for any given $\bz$, consider
the one-to-one mapping $t\mapsto v(\bz,t)$. Recall that by
Assumption \ref{assum:h} and the inverse function theorem,
$v(\bz,t)$ is continuously differentiable in $t$, and the Jacobian
\[\left|{d v(\bz,t)\over dt}\right|={1\over |\partial_m
h(\bz,v(\bz,t))|}.\]Then by the change-of-variables formula, it can
be easily verified that
\begin{eqnarray}
\int_{c_0}^{c_1}{f(\bz,v(\bz,t))\over \tilde J(\bz,v(\bz,t))}\,dt =
\int_{{\cal V}(\bz)}f(\bz,s)\,ds,\label{eqn:Case1De}
\end{eqnarray}
where ${\cal V}(\bz)$ denotes the set in which $x_m$ takes values,
i.e., ${\cal V}(\bz)=\{s: (\bz,s)\in \Omega\}$. From
(\ref{eqn:Case1De}), it can be seen that the denominator term of
$\widetilde \nu$ is indeed the density function of $\bX_{-m}$
evaluated at $\bz$. It can be evaluated based on the integral on the
RHS of (\ref{eqn:Case1De}), which can be efficiently done by
one-dimensional integration methods and does not require knowing the
function $v(\bz,t)$.

We close the discussion for Case 1 by a remark on Condition
\ref{con:C1}. Many functions of $h$ in practical applications may
satisfy Condition \ref{con:C1}, e.g., $h(\bx)=x_1+x_2$ and
$h(\bx)=x_1^2 + \exp(x_1^3x_2)$. In practice, the requirement of the
mapping being one-to-one in Condition \ref{con:C1} may be violated
for some functions, and some modifications may fix this issue. For
instance, the mapping is not one-to-one when $h(\bx)=x_1+x_2^2$. In
this case, we may either change the mapping to $\bx\mapsto
(x_2,h(\bx))$, or divide the support of $\bX$ into two parts,
$\Omega\cap \{\bx: x_2\ge 0\}$ and $\Omega\cap \{\bx: x_2<0\}$, and
apply Proposition \ref{pro:Simple} to each part separately.

However, such modifications may not work for some cases where
Condition \ref{con:C1} is violated. For instance, in a commonly
encountered case where $h(\bx)=\max(x_1,\dots,x_m)$, the mapping in
Condition \ref{con:C1} is obviously not one-to-one. To deal with
such cases, we study another sufficient condition in the following
subsection.

\subsection{Case 2}\label{sec:Case2}

We consider a function of $h$ that is homogeneous, i.e.,
\begin{condition}\label{con:C2}
The function $h$ is homogeneous, i.e., $h(t\bx)=t\bx$ for $t>0$.
\end{condition}

Homogeneous functions include, for example, the maximum function
$h(\bx)=\max(x_1,\dots,x_m)$, the minimum function
$h(\bx)=\min(x_1,\dots,x_m)$, linear functions
$h(\bx)=a_1x_1+\dots+a_mx_m$. These homogeneous functions find
important applications in a wide range of areas in operations
research.

For homogeneous functions of $h$, we consider the following mapping
\begin{equation}\label{eqn:MappingBar}
\bar u: \bx\mapsto
\left(x_1/h(\bx),\dots,x_m/h(\bx),h(\bx)\right),\quad \bar \varphi:
\by\mapsto(y_1,\dots,y_{m}),
\end{equation}
for $\by = (y_1,\dots,y_m,y_{m+1})$.

In this case, we work with a modified version of $\bX$'s support,
$\Omega'=\Omega\backslash\{\bx\in\Omega: h(\bx)=0\}$, to avoid zeros
values for denominators. When $h$ is continuously differentiable,
this modification does not change the result of our analysis. This
is because $\bX$ is a continuous random vector, and thus
$\Pr(h(\bX)=0)=0$. Therefore, ignoring the set $\{\bx\in\Omega:
h(\bx)=0\}$ does not affect our analysis on integrations.

The key condition of Theorem \ref{thm:Sensitivity} that ${\cal
H}^1\left(\bar\varphi^{-1}\{\bz\}\right)>0$ for any given $\bz$ is
satisfied when $h$ is homogeneous and the mappings are set as in
(\ref{eqn:MappingBar}). To see this, we only need to verify the
condition for two possible scenarios. The first scenario considers a
given $\bz$ for which there exists an $\bx_0$ such that
$\bz=\bx_0/h(\bx_0)$ and $h(\bx_0)>0$. Since $(c_0,c_1)$ is the
interval in which $h(\bx)$ takes values, we can see that $c_1>0$,
and
\[\bar\varphi^{-1}\{\bz\}=\{(\bz,t)\in \bar
u\{\Omega'\}:t\in (c_0,c_1)\}=\{(\bz,t):t\in (c_0,c_1),
h(t\bz)=t\}=\{(\bz,t):t\in (0,c_1)\},\]where the last equality
follows from Condition \ref{con:C2} that
$h(t\bz)=th(\bz)=th(\bx_0)/h(\bx_0)=t$ for all $t>0$. Therefore,
$\bar\varphi^{-1}\{\bz\}$ is a straight line and has positive
length, i.e, ${\cal H}^1(\bar\varphi^{-1}\{\bz\})>0$. The second
scenario considers a given $\bz$ for which there exists an $\bx_0$
such that $\bz=\bx_0/h(\bx_0)$ and $h(\bx_0)<0$. Then $c_0<0$, and
\[\bar\varphi^{-1}\{\bz\}=\{(\bz,t)\in \bar
u\{\Omega'\}:t\in (c_0,c_1)\}=\{(\bz,t):t\in (c_0,c_1),
h(t\bz)=t\}=\{(\bz,t):t\in (c_0,0)\},\]where the last equality
follows from Condition \ref{con:C2} that for any $t<0$,
\[h(t\bz)=-th(-\bz)=-th(\bx_0/(-h(\bx_0)))=-th(\bx_0)/(-h(\bx_0))=t,\]
because $h(\bx_0)<0$. Therefore, $\bar\varphi^{-1}\{\bz\}$ is a
straight line and has positive length, i.e., ${\cal
H}^1(\bar\varphi^{-1}\{\bz\})>0$.

It can be easily check that the mapping $\bar u$ is continuously
differentiable because $h$ is continuously differentiable. It is
also one-to-one, and $\bar u^{-1}(\by)=t\bz$ for any $\by=(\bz,t)\in
\bar u\{\Omega'\}$. Moreover, the Jacobian of $\bar u$ has a neat
form. By elementary algebra, the $m$-dimensional Jacobian of $\bar
u$ is
\[\bar J_m^{\bar u}(\bx) = \sqrt{(\partial_{1} h(\bx))^2+\dots+(\partial_{m}h(\bx))^2\over h^{2(m-1)}(\bx)},\]where $\partial_i h$ denotes the partial derivative of
$h(\bx)$ w.r.t. $x_i$ for $i=1,\dots,m$. Define
\[\bar J_m(\bz,t) = {1\over |t|^{m-1}}\sqrt{\left(\partial_{1}
h(t\bz)\right)^2+\dots+\left(\partial_{m}h(t\bz)\right)^2}.\]Then it
can be verified that
\[\bar J_m^{\bar u}(\bar u^{-1}(\by))=\bar
J_m(\bz,t)\left/|h(\text{sign}(t)\bz)|^{m-1}\right., \quad \forall\
\by=(\bz,t)\in u\{\Omega'\},\]where $\text{sign}(t)$ is a sign
function that is equal to $1$ if $t>0$ and $-1$ if $t<0$.

Applying the result of Theorem \ref{thm:Sensitivity}, we arrive at
an expression of $\alpha'(\xi)$, which is summarized in the
following proposition, whose proof is a direct application of
Theorem \ref{thm:Sensitivity} and thus omitted.

\begin{pro}\label{pro:main}
If Assumptions \ref{assum:Density}-\ref{assum:h} and Condition
\ref{con:C2} are satisfied, then
\[\alpha(\xi)=
\E\left[\bar w\left(\bX/h(\bX);\xi\right)\right],\]where
\[\bar w(\bz;\xi)=\left.\int_{c_0}^{\xi}|h(\text{sign}(t)\bz)|^{m-1}g(t\bz)f(t\bz)/\bar
J_m(\bz,t)dt\right/
\int_{c_0}^{c_1}|h(\text{sign}(t)\bz)|^{m-1}f(t\bz)/\bar
J_m(\bz,t)dt.\]If, in addition, $\bar
w\left(\bX/h(\bX);\cdot\right)$ satisfies a Lipschitz continuity
condition w.r.t. $\xi$ as in Theorem \ref{thm:Sensitivity}, then
\[\alpha'(\xi) =
\E\left[\bar\nu(\bX/h(\bX);\xi)\right],\]where
\[\bar\nu(\bz;\xi)=\left.{|h(\text{sign}(\xi)\bz)|^{m-1}g(\xi\bz)f(\xi\bz)\over
\bar
J_m(\bz,\xi)}\right/\int_{c_0}^{c_1}{|h(\text{sign}(t)\bz)|^{m-1}f(t\bz)\over
\bar J_m(\bz,t)}dt.\]
\end{pro}

%It should be emphasized that the result of Proposition
%\ref{pro:main} is independent of the structure of the simulation
%model of $\bX$. This property stands out in sharp contrast to
%conventional CMC which is highly problem dependent when deriving
%formulas the conditional expectation.

To illustrate how Proposition \ref{pro:main} can be applied to
derive estimators for different functions of $h$. We consider three
examples, where $h$ is chosen to be a maximum function, a linear
function, and a quadratic function, respectively, i.e.,
$h(\bx)=\max(x_1,\ldots,x_m)$, $h(\bx)=x_1+\ldots+x_m$ and
$h(\bx)=x_1^2+\dots+x_m^2$.

\begin{exmp}\label{example:Sensitivity}

\textnormal{When $h(\bx)=\max(x_1,\dots,x_m)$,
$\Omega=\mathbb{R}^m_+$ and $\xi>0$, it can be easily verified that
$\bar J_m(\bz,t)=1/t^{m-1}$. Applying Proposition \ref{pro:main}, we
have
\begin{eqnarray*}
\bar w(\bz;\xi) =
\left.\int_{0}^{\xi}g(t\bz)t^{m-1}f(t\bz)dt\right/\int_{0}^{\infty}t^{m-1}f(t\bz)dt,
\end{eqnarray*}
and
\begin{eqnarray*}
\bar
\nu(\bz;\xi)=g(\xi\bz)\xi^{m-1}f(\xi\bz)\left/\int_{0}^{\infty}t^{m-1}f(t\bz)dt\right..
\end{eqnarray*}}

\textnormal{When $h(\bx)=x_1+\dots+x_m$ and $\Omega=\mathbb{R}^m_+$,
it can be verified that $\bar J_m(\bz,t)=\sqrt{m}/t^{m-1}$, and
$\bar w(\bz;\xi)$ and $\bar \nu(\bz;\xi)$ have the same forms as
those for the maximum function.}

\textnormal{When $h(x)=x_1^2+\ldots+x_m^2$ and
$\Omega=\mathbb{R}^m$, we work with $\sqrt{h(\bx)}$ which satisfies
Condition \ref{con:C2}. Then applying Proposition \ref{pro:main}
leads to
\begin{eqnarray*}
\bar w(\bz;\xi) =
\left.\int_{0}^{\xi^2}g(t\bz)t^{m-1}f(t\bz)dt\right/\int_{0}^{\infty}t^{m-1}f(t\bz)dt,
\end{eqnarray*}
and
\begin{eqnarray*}
\bar \nu(\bz;\xi)=2\xi
g(\xi^2\bz)\xi^{2m-2}f(\xi^2\bz)\left/\int_{0}^{\infty}t^{m-1}f(t\bz)dt\right..
\end{eqnarray*}}
\end{exmp}

\section{Applications}\label{sec:Application}

\subsection{Estimating Price Sensitivities for Financial Options}\label{Greeks}

In financial risk management, the Greek letters of options play an
important role in constructing hedging strategies. Mathematically,
the Greek letters of an option are defined as the sensitivities of
the option price with respect to market parameters such as
underlying asset prices, volatilities and risk-free interest rate.
When the payoff of the option is discontinuous, estimating the Greek
letters has been a challenging problem in simulation. In what
follows, we apply the change-of-variables approach to estimate the
Greek letters for options with discontinuous payoffs.

Typically, the (discounted) payoff of an option with a discontinuous
payoff is of the form
\begin{equation}
l(\bX)\prod_{i=1}^{q}1_{\{h_i(\bX)\le a_i\}},\label{eqn:Payoff}
\end{equation}
for Lipschitz continuous functions $l$ and $h_i$'s, where
$\bX=(X_1,\dots,X_m)$ is a random vector that represents the price
dynamics of the underlying asset and depends on a market parameter
$\theta$, and for notational ease we suppress this dependence when
there is no confusion. Without loss of generality, we assume that
$\theta$ is a scalar. Based on Theorem 1 of Liu and Hong (2011), a
Greek letter associated with the market parameter $\theta$ is
represented as
\begin{eqnarray*}
\lefteqn{{d\over
d\theta}\E\left[l(\bX)\prod_{i=1}^{q}1_{\{h_i(\bX)\le
a_i\}}\right]}\\
&=&\E\left[\partial_\theta l(\bX)\prod_{i=1}^{q}1_{\{h_i(\bX)\le
a_i\}}\right]-\sum_{i=1}^{q}{d\over d
a_i}\E\left[l(\bX)\partial_\theta h_i(\bX)\prod_{k\neq
i}1_{\{h_k(\bX)\le a_k\}} \cdot 1_{\{h_i(\bX)\le a_i\}}\right],
\end{eqnarray*}
where $\partial_\theta$ denotes the operator of takin derivative
with respect to $\theta$.

Note that in a simulation run, both $\{l(\bX),h_i(\bX),1\le i\le
q\}$ and their pathwise derivatives $\{\partial_\theta
l(\bX),\partial_\theta h_i(\bX),1\le i\le q\}$ are usually readily
computable; see, e.g., Broadie and Glasserman (1996). The first term
on the RHS of the above equation can then be easily estimated by a
sample-mean estimator. In the rest of this subsection, we will
discuss how to estimate the second term using the
change-of-variables approach.

To simplify notation, we let
\[g_i(\bX)=l(\bX)\partial_\theta
h_i(\bX)\prod_{k\neq i}1_{\{h_k(\bX)\le a_k\}}.\]Then the problem is
reduced to how to estimate
\[\beta_i \triangleq {d\over d a_i}\E\left[g_i(\bX)\cdot 1_{\{h_i(\bX)\le
a_i\}}\right],\quad i = 1,\dots,q.\]

For many options traded in financial markets, the function
$h_i(\bX)$ in their payoffs is often in the form of $h_i(\bX)=X_m$,
$h_i(\bX)=\sum_{i=1}^{m}X_i/m$, $h_i(\bX)=\max(X_1,\dots,X_m)$, or
$h_i(\bX)=\min(X_1,\dots,X_m)$; see, e.g., Tong and Liu (2016) for
more detailed discussions. It can be easily seen that these are all
homogenous functions and satisfy Condition \ref{con:C2}. Applying
Proposition \ref{pro:main}, we have
\begin{eqnarray*}
\beta_i={d\over d a_i}\E\left[g_i(\bX)\cdot 1_{\{h_i(\bX)\le
a_i\}}\right]= \E\left[\nu\left(\bZ;a_i\right)\right],
\end{eqnarray*}
where $\bZ\triangleq \left(X_1/h_i(\bX),\dots,X_m/ h_i(\bX)\right)$,
and
\[\nu(\bz;a_i)=g_i(\bz a_i)a_i^{m-1}f(\bz
a_i)\left/\int_{0}^{\infty}t^{m-1}f(t\bz)dt\right..\]It turns out
that the same form of $\nu$ applies to all $\beta_i$, $i=1,\dots,q$.

Deriving an explicit formula of $\nu(\bz;a_i)$ requires knowing the
joint density function $f$, which is available for many commonly
used pricing models, for instance, when $\bX$ represents the
underlying asset prices observed at different time points under the
Black-Scholes model. As a remark, it is worth pointing out that this
requirement can be further relaxed. Indeed, one may replace $f$ by a
conditional density of $\bX$ given some variables. This relaxation
offers considerable flexibility during implementation, as
conditional densities can often be obtained for most, if not all,
pricing models for financial options. As an illustrative example, we
derive $\nu(\bz;a_i)$ for the more complex variance gamma model
using a conditional density where an explicit joint density is not
available; see Section \ref{sec:VG} of the online supplement for
more details.

When estimating $\beta_i$, we generate $n$ identically and
independently distributed (i.i.d.) observations of $\bX$, denoted by
$\{\bX_1,\dots,\bX_n\}$, and compute
\[\bZ_k=\left({X_{k,1}\over
h_i(\bX_k)},\dots,{X_{k,m}\over h_i(\bX_k)}\right),\quad
k=1,\dots,n,\]where $\bX_k=(X_{k,1},\dots,X_{k,m})$.

Then $\beta_i$ can be estimated by
\[\bar M_n = {1\over n}\sum_{k=1}^{n}\nu(\bZ_k;a_i).\]
It is worth mentioning that $\bar M_n$ is an unbiased estimator,
thus it has desirable asymptotic properties as a typical sample-mean
estimator. Although it involves another random vector $\bZ$, the
estimator does not require any change of probability measures in the
simulation, because $\bZ$ is readily computable once $\bX$ is
generated.

\subsubsection{Numerical Experiments}\label{GreekExperiments}

We consider two pricing models, the Black-Scholes (BS) model and the
variance gamma (VG) model, to conduct numerical experiments and
examine the performances of the estimators. In particular, the price
of the underlying asset is monitored at $m$ time points
$\{t_1<\dots<t_m\}$ evenly spaced over $(0,T)$, i.e., $t_i=iT/m$ for
$i=1,\ldots,m$, where $T$ is the maturity date of the option. To
simplify notation, we let $X_i$ denote the underlying asset price at
$t_i$ for $i=0,\dots,m$.

Under the BS model, the price of the underlying asset is governed by
a geometric Brownian motion, i.e.,
\[X_{i+1}=X_i\exp\left((r-\sigma^2/2)T/m+\sigma \sqrt{T/m}N_{i+1}\right),\quad i=0,\dots,m-1,\]where
$\{N_1,\dots,N_m\}$ are independent standard normal random
variables, $r$ and $\sigma$ denote the risk-free interest rate and
volatility of the underlying asset. Denote the initial price of the
underlying asset by $X_0=x_0$. The joint density of $\bX$ is then
\[f(\bx)=\prod_{i=0}^{m-1}{1\over \sigma\sqrt{T/m}x_{i+1}}\phi\left({1\over
\sigma\sqrt{T/m}}\left(\log(x_{i+1}/x_i)-(r-\sigma^2/2)T/m\right)\right),\]where
$\phi(\cdot)$ denotes the standard normal density function.
%\[f(x_1,\dots,x_m)=\left({1\over \sqrt{2\pi T/m}\sigma}\right)^m{1\over x_1\dots x_m}\exp\left(-{1\over 2\sigma^2T/m}\sum_{i=1}^{m}\left(\log (x_i/x_{i-1})-(r-\sigma^2/2)T/m\right)^2\right).\]

Compared to the BS model, the VG model is a pure jump process and
allows for more flexible skewness and kurtosis; see, e.g., Madan et
al. (1998). A discretization of the VG model is given by (Fu 2000),
\begin{equation*}
X_{i+1}=X_i\exp(\mu T/m +\theta G_{i+1}+\sigma
\sqrt{G_{i+1}}N_{i+1}), \quad i = 0,\dots,m-1,
\end{equation*}
where $\{G_1,\dots,G_m\}$ are independent gamma random variables
with scale parameter $T/(m\beta)$ and shape parameter $\beta$,
$\{N_1,\dots,N_m\}$ are independent standard normal random
variables, $\mu=r + 1/\beta\log (1-\theta\beta-\sigma^2\beta/2)$,
and $\theta$ and $\sigma$ are parameters of the model. Then the
conditional density of $\bX$ given $\bG\triangleq (G_1,\dots,G_m)$
is
\[f(\bx|\bG)=\prod_{i=0}^{m-1}{1\over x_{i+1}\sigma\sqrt{G_{i+1}}}\phi\left({1\over \sigma\sqrt{G_{i+1}}}\left(\log(x_{i+1}/x_i)-\mu T/m-\theta G_{i+1}\right)\right).\]

Under each of the BS and VG models, we consider three options with
discontinuous payoffs, including a digital option with discounted
payoff $e^{-rT}1_{\{X_m\ge K\}}$, an Asian digital option with
discounted payoff $e^{-rT}1_{\{\sum_{i=1}^{m}X_i/m\ge K\}}$, and a
barrier call option with discounted payoff
$e^{-rT}(X_m-K)^+1_{\{\max(X_1,\dots,X_m)\le \kappa\}}$, where $K$
and $\kappa$ denote the strike price and the barrier respectively.
For each of the above options, we estimate the Greek letters {\tt
delta} and {\tt gamma}, i.e., the first- and second-order
derivatives of the option price w.r.t. $x_0$, and {\tt theta} and
{\tt vega}, i.e., the first-order derivatives of the option price
w.r.t. $T$ and $\sigma$, respectively. We compare the proposed
change-of-variables estimators to existing ones in the literature,
including the likelihood ratio method and conventional CMC when
applicable. Detailed derivation of various estimators is provided in
Section \ref{sec:OnlineDerive} of the online supplement.

In all experiments, we set the sample size as $n=10^5$. To examine
the performance of an estimator $\bar M_n$, we use its relative
error, defined as the ratio of the standard deviation of $\bar M_n$
to the absolute value of the quantity being estimated. True values
of the quantities being estimated are either computed by closed-form
formulas when available, or approximated using other existing
methods with an extremely large sample size ($10^9$).

For the options under the BS model, we let $x_0=K=100$,
$\kappa=120$, $r=5\%$, $\sigma=0.3$, and $T=1$. We vary the number
of discretization steps, $m$, to examine its impact on the
performances of various estimators, including the likelihood ratio
(LR) estimator, conventional CMC estimator, and the proposed
change-of-variables (CoV) estimator. Comparison results for the
digital, Asian and barrier options are summarized in Tables
\ref{tab:BS_digital_greeks}, \ref{tab:BS_asian_greeks} and
\ref{Table:BS_barrier_greeks}, respectively. From these tables it
can be seen that the proposed CoV estimators have the best
performances in all settings. Its improvement upon existing methods
can be dramatic. For instance, when estimating {\tt gamma} for the
Asian digital option with $m=100$, relative errors of the LR and
conventional CMC estimators are over $10$ and $1000$ times of that
of the CoV estimator respectively, implying that sample sizes of the
LR and CMC estimators have to be as large as $100$ and $10^6$ times
of that of the CoV estimator in order to achieve the same level of
accuracy. For the barrier option, it is not clear how conventional
CMC estimators can be derived, while the proposed CoV estimators
perform very well.

\begin{table}[h!]
\vspace*{-1pt}
\begin{center}
\caption{Relative errors (\%) of various Greek estimators for the
digital option under the BS model} \label{tab:BS_digital_greeks}
\def\temptablewidth{1\textwidth}
{\rule{\temptablewidth}{1pt}}
\begin{tabular*}{\temptablewidth}{@{\extracolsep{\fill}}rcccccccccccccccc}
% {\rule[-1mm]{0mm}{5mm}}
% & \multicolumn{15}{c}{Relative Error}\\
% \cline{2-16}{\rule[-1mm]{0mm}{7mm}}
 & \multicolumn{3}{c}{\tt delta} & & \multicolumn{3}{c}{\tt vega} & & \multicolumn{3}{c}{\tt theta} & & \multicolumn{3}{c}{\tt gamma}\\
 \cline{2-4}
 \cline{6-8}
 \cline{10-12}
 \cline{14-16}
 {\rule[-1mm]{0mm}{7mm}}
 $m$ & LR & CMC & CoV &  & LR & CMC & CoV &  & LR & CMC & CoV & &  LR & CMC & CoV\\
 \cline{2-16}{\rule[-1mm]{0mm}{7mm}}
$10$ & 1.8 & 0.4 &0.4&           & 7.9& 0.4& 0.4&         & 23.0 &0.4 &0.4&         &24.9 & 3.3& 3.3\\
$50$ &4.0 &0.6& 0.6&              & 17.8& 0.6& 0.6&        & 51.1& 0.4& 0.4&         & 125& 11.1& 11.1 \\
$100$ & 5.6 & 0.8& 0.8&           & 25.0& 0.8& 0.8& & 71.8& 0.4&
0.4&       & 252& 18.7& 18.7
\end{tabular*}
{\rule{\temptablewidth}{1pt}}
\end{center}
\end{table}

\begin{table}[h!]
\vspace*{-1pt}
\begin{center}
\caption{Relative errors (\%) of various Greek estimators for the
Asian digital option under the BS model} \label{tab:BS_asian_greeks}
\def\temptablewidth{1\textwidth}
{\rule{\temptablewidth}{1pt}}
\begin{tabular*}{\temptablewidth}{@{\extracolsep{\fill}}rcccccccccccccccc}
% {\rule[-1mm]{0mm}{5mm}}
% & \multicolumn{15}{c}{Relative Error}\\
% \cline{2-16}{\rule[-1mm]{0mm}{7mm}}
 & \multicolumn{3}{c}{\tt delta} & & \multicolumn{3}{c}{\tt vega} & & \multicolumn{3}{c}{\tt theta} & & \multicolumn{3}{c}{\tt gamma}\\
 \cline{2-4}
 \cline{6-8}
 \cline{10-12}
 \cline{14-16}
 {\rule[-1mm]{0mm}{7mm}}
 $m$ & LR & CMC & CoV &          & LR & CMC & CoV &           & LR & CMC & CoV &          &  LR & CMC & CoV\\
 \cline{2-16}{\rule[-1mm]{0mm}{7mm}}
$10$ & 1.1 & 1.1 &0.2 &          & 10.0& 1.1 & 0.2 &      & 28.8 &0.7 &0.5 &         &14.9 & 86.2& 2.5\\
$50$ &2.3 &3.7& 0.4&             & 22.8& 3.7 & 0.4 &       & 64.0 & 2.1 & 0.5 &       & 74.1& 3210 & 8.6 \\
$100$ & 3.3 & 6.3 & 0.6 &       & 32.7 & 6.2 & 0.6 & & 91.3 & 3.5 &
0.6 &       & 151 & 15164 & 14.4
\end{tabular*}
{\rule{\temptablewidth}{1pt}}
\end{center}
\end{table}

\begin{table}[h!]
\vspace*{-1pt}
\begin{center}
\caption{Relative errors (\%) of various Greek estimators for the
barrier call option under the BS model}
\label{Table:BS_barrier_greeks}
\def\temptablewidth{0.8\textwidth}
{\rule{\temptablewidth}{1pt}}
\begin{tabular*}{\temptablewidth}{@{\extracolsep{\fill}}rcccccccccccc}
% {\rule[-1mm]{0mm}{5mm}}
% & \multicolumn{11}{c}{Relative Error}\\
% \cline{2-12}{\rule[-1mm]{0mm}{7mm}}
 & \multicolumn{2}{c}{\tt delta} & & \multicolumn{2}{c}{\tt vega} & & \multicolumn{2}{c}{\tt theta} & & \multicolumn{2}{c}{\tt gamma}\\
 \cline{2-3}
 \cline{5-6}
 \cline{8-9}
 \cline{11-12}
 {\rule[-1mm]{0mm}{7mm}}
 $m$ & LR &  CoV &  & LR &  CoV &  & LR &  CoV & &  LR &  CoV\\
 \cline{2-12}{\rule[-1mm]{0mm}{7mm}}
$10$ & 5.5  &5.6 & & 2.6 &0.5 & & 2.6  &0.7 & &5.6 & 3.0\\
$50$ &11.1 &5.6 & & 5.9 &0.8 & &  5.9 & 1.1& & 41.7 & 11.2 \\
$100$ & 15.1 & 6.2 & & 8.3 & 0.9 & &  8.3 & 1.3 & & 84.9 & 19.8
\end{tabular*}
{\rule{\temptablewidth}{1pt}}
\end{center}
\end{table}

For options under the VG model, we let $x_0 = K =100$, $\kappa =
120$, $r =5\%$, $\sigma = 0.2$, $\beta = 10$, $\theta = -0.2$, and
$T=1$. Comparison results for the digital, Asian and barrier options
are summarized in Tables \ref{table:VG_digital_greeks},
\ref{table:VG_asian_greeks} and \ref{table:VG_barrier_greeks},
respectively. From the tables it can be seen that the proposed CoV
estimators significantly outperform the existing estimators in many
cases. For instance, when estimating {\tt theta} for the Asian
option with $m=100$, relative errors of the LR and conventional CMC
estimators are over $280$ and $10$ times of that of the CoV
estimator.

\begin{table}[h!]
\vspace*{-1pt}
\begin{center}
\caption{Relative errors (\%) of various Greek estimators for the
digital option under the VG model} \label{table:VG_digital_greeks}
\def\temptablewidth{1\textwidth}
{\rule{\temptablewidth}{1pt}}
\begin{tabular*}{\temptablewidth}{@{\extracolsep{\fill}}rcccccccccccccccc}
% {\rule[-1mm]{0mm}{5mm}}
% & \multicolumn{15}{c}{Relative Error}\\
% \cline{2-16}{\rule[-1mm]{0mm}{7mm}}
 & \multicolumn{3}{c}{\tt delta} & & \multicolumn{3}{c}{\tt vega} & & \multicolumn{3}{c}{\tt theta} & & \multicolumn{3}{c}{\tt gamma}\\
 \cline{2-4}
 \cline{6-8}
 \cline{10-12}
 \cline{14-16}
 {\rule[-1mm]{0mm}{7mm}}
 $k$ & LR & CMC & CoV &  & LR & CMC & CoV &  & LR & CMC & CoV & &  LR & CMC & CoV\\
 \cline{2-16}{\rule[-1mm]{0mm}{7mm}}
$10$ & 1.7 &0.4 &0.4 &  & 14.9& 0.5&0.5& &31.1 &0.3 &0.3&   &383 &52.9& 52.8\\
$50$ &3.9 &0.7& 0.7&    &32.1&0.7& 0.7&   &68.5&0.5& 0.5&  &1343&121& 121 \\
$100$ &5.4 &0.8 &0.8&  &44.8&0.8& 0.8& &97.2&0.6& 0.6& &2579&197&
197
\end{tabular*}
{\rule{\temptablewidth}{1pt}}
\end{center}
\end{table}

\begin{table}[h]
\vspace*{-1pt}
\begin{center}
\caption{Relative errors (\%) of various Greek estimators for the
Asian digital option under the VG model}
\label{table:VG_asian_greeks}
\def\temptablewidth{1\textwidth}
{\rule{\temptablewidth}{1pt}}
\begin{tabular*}{\temptablewidth}{@{\extracolsep{\fill}}rcccccccccccccccc}
% {\rule[-1mm]{0mm}{5mm}}
% & \multicolumn{15}{c}{Relative Error}\\
% \cline{2-16}{\rule[-1mm]{0mm}{7mm}}
 & \multicolumn{3}{c}{\tt delta} & & \multicolumn{3}{c}{\tt vega} & & \multicolumn{3}{c}{\tt theta} & & \multicolumn{3}{c}{\tt gamma}\\
 \cline{2-4}
 \cline{6-8}
 \cline{10-12}
 \cline{14-16}
 {\rule[-1mm]{0mm}{7mm}}
 $k$ & LR & CMC & CoV &  & LR & CMC & CoV &  & LR & CMC & CoV & &  LR & CMC & CoV\\
 \cline{2-16}{\rule[-1mm]{0mm}{7mm}}
$10$ & 1.0 &1.1 &0.2 &  & 14.9& 1.4&0.3& &35.8 &0.9 &0.2&   &34.8 &188& 5.6\\
$50$ &2.3 &3.9& 0.5&    &33.2&4.2& 0.5&   &80.0&2.5& 0.3&  &171&7007& 19.8 \\
$100$ &3.2 &6.7&0.6&  &46.8&6.8& 0.6&   &112&4.2& 0.4& &323&34262&
32
\end{tabular*}
{\rule{\temptablewidth}{1pt}}
\end{center}
\end{table}

\begin{table}[h!]
\vspace*{-1pt}
\begin{center}
\caption{Relative errors (\%) of various Greek estimators for the
barrier call option under the VG model}
\label{table:VG_barrier_greeks}
\def\temptablewidth{0.8\textwidth}
{\rule{\temptablewidth}{1pt}}
\begin{tabular*}{\temptablewidth}{@{\extracolsep{\fill}}rcccccccccccc}
% {\rule[-1mm]{0mm}{5mm}}
% & \multicolumn{11}{c}{Relative Error}\\
% \cline{2-12}{\rule[-1mm]{0mm}{7mm}}
 & \multicolumn{2}{c}{\tt delta} & & \multicolumn{2}{c}{\tt vega} & & \multicolumn{2}{c}{\tt theta} & & \multicolumn{2}{c}{\tt gamma}\\
 \cline{2-3}
 \cline{5-6}
 \cline{8-9}
 \cline{11-12}
 {\rule[-1mm]{0mm}{7mm}}
 $k$ & LR &  CoV &  & LR &  CoV &  & LR &  CoV & &  LR &  CoV\\
 \cline{2-12}{\rule[-1mm]{0mm}{7mm}}
$10$ & 8.7  &5.6 & & 3.1 &  1.1 & & 11.0  &0.6 & &7.3 & 2.3\\
$50$ &50.4 &19.5 & & 6.1 &1.6& &  22.3 & 0.9& & 40.0 & 7.8 \\
$100$ & 106 & 35.2& & 8.3 & 1.9 & &  29.7 & 1.0 & & 72.5 & 12.9
\end{tabular*}
{\rule{\temptablewidth}{1pt}}
\end{center}
\end{table}

\subsection{Estimating Gradient of Chance Constrained Programs}

Consider a chance constrained program
\begin{eqnarray}\label{CCP}
&\displaystyle\mathop{\mbox{minimize}}_{\bt} & r(\bt)\\
&\mbox{ subject to} & \Pr\{L(\bt,\bX)\le 0\}\ge \beta,\nonumber
\end{eqnarray}
where $\bt\triangleq (t_1,\ldots,t_m)^T\in \mathbb{R}^m$ represents
the vector of decision variables, $r(\bt)$ is a deterministic
objective function, $\bX=(X_1,\ldots,X_m)^T$ is an $m$-dimensional
random vector, $L$ is a known function, and $\beta\in (0,1)$ is a
parameter specified by the modeler.

Sample average approximation (SAA) is a popular method for solving
the chance constrained program (\ref{CCP}); see, e.g., Luedtke and
Ahmed (2008) and Pagnoncelli et al. (2009). To use SAA, it would be
helpful if one has an estimate of the gradient of the constraint,
i.e.,
\[\nabla_{\bt} \Pr\{L(\bt,\bX)\le 0\}.\]To illustrate how the
change-of-variables approach can be applied to estimate this
gradient, we consider an example where $L$ is linear, i.e.,
\[L(\bt,\bX)=\bt^T\bX-b,\]where $b$ is a given constant.

By using Theorem 1 of Liu and Hong (2011), it can be shown that for
$i=1,\ldots,m$,
\[{\partial\over \partial t_i}\Pr\{L(\bt,\bX)\le 0\}=-{\partial\over \partial b} \E\left[X_i\cdot 1_{\left\{\bt^T\bX\le b\right\}}\right].\]
Let $h(\bx)=\bt^T\bx$. It is straightforward that Condition
\ref{con:C2} is satisfied. Applying Proposition \ref{pro:main}, we
have
\begin{eqnarray}
\nabla_{\bt}\Pr\left\{\bt^T\bX\le
b\right\}=\E\left[\nu(\bZ;b)\right], \label{eqn:EstCCP}
\end{eqnarray}
where
\[\nu(\bZ;b)=-{b\bZ f(b\bZ)|b|^{m-1}\over \int_{-\infty}^{\infty}f(t\bZ)|t|^{m-1}dt},\]where $f$ is the density function of $\bX$,
and $\bZ=\bX/t^T\bX$.

When estimating the gradient, we generate $n$ i.i.d. observations of
$\bX$, denoted by $\{\bX_1,\dots,\bX_n\}$, and compute
$\bZ_k=\bX_k/\bt^T\bX_k$ for $k=1,\dots,n$. Then the gradient can be
estimated by
\[\bar M_n = {1\over n}\sum_{k=1}^{n}\nu(\bZ_k;b).\]

\subsubsection{Numerical Experiments}

Numerical experiments are conducted to illustrate the performance of
the gradient estimator resulting from (\ref{eqn:EstCCP}). In
particular, we estimate $\partial_{t_1}\Pr\left\{\bt^T\bX\le
b\right\}$ at $\bt=(1,1,\ldots,1)^T$. We vary $b$ to examine the
performance of the estimator in different settings. More
specifically, we choose $b$ such that the probability of
$\{\bt^T\bX\le b\}$ takes values in $\{0.90,0.95,0.99\}$.

Consider two cases, where $\bX$ follows a multivariate normal and a
multivariate Student's t-distribution respectively. Detailed
derivation of explicit formulas for $\nu$ in (\ref{eqn:EstCCP}) is
provided in Section \ref{sec:OnlineCCP} of the online supplement. In
our numerical setting, we let the means of the distributions be zero
and the covariance matrix $\Sigma$ is specified in a way that
individual $X_i$ has a unit variance while every pair $(X_i,X_j)$
has a correlation $\rho$ when $i\neq j$. We let $\rho=0.3$ in the
numerical experiments, and the degrees of freedom be $4$ for the
multivariate t-distribution.

We compare the estimators resulting from the proposed
change-of-variables (CoV) approach to conventional CMC estimators.
To measure the performance of an estimator, we report its relative
error. In addition, we report the ratio of the relative error of the
conventional CMC estimator to that of the CoV estimator. In all
experiments, we set the sample size $n$ as $10^5$. Numerical results
are summarized in Tables \ref{Table:NormalCCP} and \ref{Table:tCCP}
for the multivariate normal and t cases, respectively. From the
tables it can be seen that the proposed CoV estimator performs
better in all settings. Its performance is significantly better than
conventional CMC in some cases. For instance, for the multivariate t
case with $m=50$ and a probability level of $99\%$, the ratio of
their relative errors is over $10$, implying that the sample size of
the conventional CMC estimator has to be as large as $100$ times of
that of the CoV estimator in order to achieve the same level of
accuracy.

\begin{table}[hhh]
\caption{Relative errors (\%) and performance ratios of the
conventional CMC and CoV estimators for the multivariate normal
case. } \vspace*{-12pt}\label{Table:NormalCCP}
\begin{center}
\def\temptablewidth{1.0\textwidth}
{\rule{\temptablewidth}{1pt}}
\begin{tabular*}{\temptablewidth}{@{\extracolsep{\fill}}lrrrrrrrrrrr}
 {\rule[-1mm]{0mm}{5mm}}
& \multicolumn{3}{c}{$m=5$} & &\multicolumn{3}{c}{$m=10$} & &\multicolumn{3}{c}{$m=50$}\\%[0.1ex]
\cline{2-4}  \cline{6-8}  \cline{10-12}\\
$\Pr\left\{\bt^T\bX\le
b\right\}$ & 0.90& 0.95 & 0.99 & & 0.90& 0.95 & 0.99 && 0.90& 0.95 & 0.99\\
\cline{2-4}  \cline{6-8}  \cline{10-12}{\rule[-1mm]{0mm}{5mm}}\\
\text{CMC}& 0.7& 0.9 & 1.7 & & 1.2 & 1.4 & 2.6 && 2.9& 3.5 & 6.4\\
\text{CoV}& 0.5& 0.5 & 0.8 & & 0.7& 0.7 & 1.0 && 1.2& 1.2 & 1.7\\
\text{Ratio}& 1.5& 1.7 & 2.3 & & 1.7& 2.0 & 2.5 && 2.5& 3.0 & 3.7
\end{tabular*} {\rule{\temptablewidth}{1pt}}
\end{center}
       \end{table}

       \begin{table}[hhh]
\caption{Relative errors (\%) and performance ratios of the
conventional CMC and CoV estimators for the multivariate t case. }
\vspace*{-12pt}\label{Table:tCCP}
\begin{center}
\def\temptablewidth{1.0\textwidth}
{\rule{\temptablewidth}{1pt}}
\begin{tabular*}{\temptablewidth}{@{\extracolsep{\fill}}lrrrrrrrrrrr}
 {\rule[-1mm]{0mm}{5mm}}
& \multicolumn{3}{c}{$m=5$} & &\multicolumn{3}{c}{$m=10$} & &\multicolumn{3}{c}{$m=50$}\\%[0.1ex]
\cline{2-4}  \cline{6-8}  \cline{10-12}\\
$\Pr\left\{\bt^T\bX\le
b\right\}$& 0.90& 0.95 & 0.99 & & 0.90& 0.95 & 0.99 && 0.90& 0.95 & 0.99\\
\cline{2-4}  \cline{6-8}  \cline{10-12}{\rule[-1mm]{0mm}{5mm}}\\
\text{CMC}& 0.9& 1.1 & 2.3 & & 1.4& 1.7 & 3.5 && 3.4& 4.2 & 8.3\\
\text{CoV}& 0.4& 0.5 & 0.5 & & 0.5& 0.6 & 0.7 && 0.6& 0.6 & 0.8\\
\text{Ratio}& 2.0& 2.4 & 4.2 & & 2.6& 3.2 & 5.2 && 5.4& 6.6 & 10.6
\end{tabular*} {\rule{\temptablewidth}{1pt}}
\end{center}
       \end{table}

\section{Concluding Remarks}\label{sec:Conclusions}

In this paper, we propose a change-of-variables approach to CMC, and
provide theoretical underpinnings. We attempt to circumvent the
difficulty on finding conditioning variables for CMC. To this end,
we provide sensitivity estimators for discontinuous integrands under
appropriate smoothness conditions. Many practical applications may
fit into our setting. We show that the proposed approach may lead to
new estimators, exemplified by two applications, including
estimation of sensitivities of financial options with discontinuous
payoffs and gradient estimation of chance constrained programs.

The change-of-variables approach might also find applications in
other simulation problems, e.g., sampling in a hyperplane in
$\mathbb{R}^n$ from a given distribution. In principle, sampling
from an properly chosen distribution in $\mathbb{R}^{n-1}$ and then
mapping the generated samples to the hyperplane might achieve the
goal. However, how to design the mapping is challenging, and it is
left as a topic for future research.

\section*{Acknowledgements}

The research of the second author is partially supported by the Hong
Kong Research Grants Council under grants CityU 155312 and 192313.

\section*{References}

\begin{hangref}

\item Asmussen, S., and P. W. Glynn. 2007. {\it Stochastic Simulation: Algorithms and
Analysis}, Springer, New York.

%\item Asmussen, S., and D. P. Kroese. 2006. Improved algorithms for
%rare event simulation with heavy tails. {\it Advances in Applied
%Probabilities}, {\bf 38}, 545-558.

\item Bernis, G., E. Gobet and A. Kohatsu-Higa. 2003. Monte Carlo
evaluation of Greeks for multidimensional barrier and lookback
options. {\it Mathematical Finance}, {\bf 13} 99-113.

\item Broadie, M., and P. Glasserman. 1996. Estimating security price derivatives using simulation. {\it Management Science}, {\bf 42} 269-285.

%\item Blanchet, J. H., J. Liu, and X. Yang. 2010. Monte Carlo for
%large credit portfolios with potentially high correlations. {\it
%Proceedings of the 2010 Winter Simulation Conference}, 2810-2820.

\item Chan, J. H. and M. Joshi. 2013. Fast Monte Carlo Greeks for financial products with discontinuous pay-offs. {\it Mathematical Finance}, {\bf 23}(3) 459-495.

\item Chen, N., and P. Glasserman. 2007. Malliavin Greeks without
Malliavin calculus. {\it Stochastic Processes and their
Applications}, {\bf 117}:1689-1723.

%\item Chan, J. H., and M. Joshi. 2015. Optimal limit methods for
%computing sensitivities of discontinuous integrals including
%triggerable derivative securities. {\it IIE Transactions}, (\bf 47),
%978-997.

\item Durrett, R. 2005. {\it Probability: Theory and Examples}, Third Edition. Duxbury Press, Belmont.

%\item Evans, L. C., and R. F. Gariepy. 2015. {\it Measure Theory and Fine Properties of Functions}, CRC Press, Boca Raton.

\item Federer, H. 1996. {\it Geometric Measure Theory}, Reprint of the 1969 Edition, Springer, Berlin.

\item Fu, M. C. 2000. Variance-Gamma and Monte Carlo.
M. C. Fu, R. A. Jarrow, J.-Y. Yen, and R. J. Elliott eds. {\it
Advances in Mathematical Finance}. Boston, Birkhauser.

%\item Fu, M. C. 2006. Gradient estimation. S. G. Henderson, B. L.
%Nelson eds. {\it Simulation: Handbooks in Operations Research and
%Management Science}. Elsevier, Amsterdam, The Netherlands.

\item Fu, M. C., L. J. Hong, and J.-Q. Hu. 2009. Conditional Monte
Carlo estimation of quantile sensitivities. {\it Management
Science}, {\bf 55}, 2019-2027.

\item Fu, M. C., and J.-Q. Hu. 1997. {\it Conditional Monte Carlo: Gradient Estimation and Optimization
Applications}, Kluwer Academic Publishers, Boston, MA.

%\item Glasserman, P. 1991. {\it Gradient Estimation via Perturbation
%Analysis}, Kluwer Academic Publishers, Boston, MA.

\item Glynn, P. W. 1987. Likelihood ratio gradient estimation: An
overview. {\it Proceedings of the 1987 Winter Simulation
Conference}, 366-374.

\item Gong, W. B., and Y. C. Ho. 1987. Smoothed perturbation
analysis of discrete-event dynamic systems. {\it IEEE Transactions
on Automatic Control}, {\bf 32}, 858-867.

\item Hammersley, J. M. 1956. Conditional Monte Carlo. {\it Journal of the
ACM}, {\bf 3} 73-76.

\item Ho, Y. C., and X.-R. Cao. 1983. Perturbation analysis and
optimization of queueing networks. {\it Journal of Optimization
Theory and Applications}, {\bf 40} 559-582.

\item Law, A. M., and W. D. Kelton. 2000. {\it Simulation Modeling and
Analysis}, Third Edition, McGraw-Hill.

\item Liu, G., and L. J. Hong. 2011. Kernel estimation of the Greeks for options with discontinuous payoffs. {\it Operations Research}, {\bf 59} 96-108.

\item L'Ecuyer, P. 1990. A unified view of the IPA, SF, and LR
gradient estimation technique. {\it Management Science}, {\bf 36},
1364-1383.

\item Luedtke, J., S. Ahmed. 2008. A sample approximation approach for
optimization with probabilistic constraints. {\it SIAM Journal on
Optimization}, {\bf 19}(2) 674-699.

\item Lyuu, Y.-D.,and H.-W. Teng. 2011. Unbiased and efficient Greeks of financial options. {\it Finance and Stochastics}, {\bf 15} 141-181.

\item Madan, D. B., P. Carr, and E. C. Chang. 1998. The variance
gamma process and option pricing. {\it European Finance Review},
{\bf 2} 79-105.

\item Mal\'y, J. 2001. Lectures on change of variables in integral.
Preprint 305, Department of Mathematics, University of Helsinki.

\item Morgan, F. 2009. {\it Geometric Measure Theory: A Beginner's
Guide}, 3rd Edition, Academic Press.

\item Pagnoncelli, B. K., S. Ahmed, A. Shapiro. 2009. Sample average
approximation method for chance constrained programming: Theory and
applications. {\it Journal of Optimization Theory and Applications},
{\bf 142} 399-416.

\item Pflug, G. C. 1988. {\it Derivatives of Probability Measures: Concepts and Applications to the Optimization of Stochastic Systems}. Springer, Berlin.

\item Pflug, G. C., and H. Weisshaupt. 2005. Probability gradient estimation by set-valued calculus and applications in network design. {\it SIAM Journal on Optimization}, {\bf15}(3) 898-914.

\item Suri, R., and M. A. Zazanis. 1988. Perturbation analysis gives
strongly consistent sensitivity estimates for the M/G/1 queue. {\it
Management Science}, {\bf 34}, 39-64.

\item Tong, S., and G. Liu. 2016. Importance sampling for option
Greeks with discontinuous payoffs. {\it INFORMS Journal on
Computing}, in press.

\item Wang Y., M. C. Fu, and S. I. Marcus. 2009. Sensitivity
analysis for barrier options. {\it Proceedings of the 2009 Winter
Simulation Conference}, 1272-1282.

\item Wang, Y., M. C. Fu, and S. I. Marcus. 2012. A new stochastic derivative estimator for discontinuous payoff functions with application to financial derivatives. {\it Operations
Research}, {\bf 60}(2) 447-460.

\item Wendel, J. G. 1957. Groups and conditional Monte Carlo. {\it The Annals of Mathematical
Statistics}, {\bf 28} 1048-1052.

\end{hangref}

%\newpage

\appendix

\section{Appendix}\label{sec:Proofs}

%\subsection{Proofs}

%\begin{figure}[h]
%   \centering
%    \includegraphics[width = 3.0in]{original.pdf}
%   \caption{Substitutability/Complementarity for Different Values of $(\mu_1, \mu_2)$  in Example \ref{example:nested} $(\mu_3 = 3)$}
%  \label{fig:nested_logit_2}
%\end{figure}
%\begin{figure}
%  \centering
%  % Requires \usepackage{graphicx}
%  \includegraphics[width=3.0in]{original.pdf}\\
%  \caption{dssf}\label{ew}
%\end{figure}

In what follows we provide the proofs of the main results. To
facilitate the proofs, an area formula and a coarea formula for
Hausdorff integrations will be used repeatedly. Interested readers
are referred to Mal\'y (2001, Theorem 4.6, Exercise 3) for details
of the area formula, and to Federer (1996, Theorems 3.2.22 and
3.2.31) for those of the coarea formula.

\begin{lemma}[Area Formula]\label{lem:areaformula}
Let $E_0 \subset \mathbb{R}^{m}$ be an open set, and $u: E_0 \mapsto
\mathbb{R}^{n}$ be a continuously differentiable one-to-one mapping,
$n\ge m$. Then for any Borel set $E\subset E_0$ we have
\begin{equation}\label{4}
  \int_{E}J_{m}^u(\bx)d\bx=\int_{u\{E\}}
  \mathcal{H}^{m}(d\by)=\mathcal{H}^{m}(u\{E\}).
\end{equation}
Furthermore,
\begin{equation}\label{4}
  \int_{E}p(u(\bx))J_{m}^u(\bx)d\bx=\int_{u\{E\}} p(\by)\mathcal{H}^{m}(d\by),
\end{equation}
for all Borel functions $p: u\{E\}\mapsto \mathbb{R}$ for which one
side exists.
\end{lemma}

\begin{lemma}[Coarea Formula]\label{lem:coareaformula}
Suppose $n\geq m\geq k\leq q$ are positive integers, $A$ is an
$(\mathcal{H}^{m},m)$ rectifiable Borel set of $\mathbb{R}^{n}$,
$E\subset \mathbb{R}^{q}$, and the mapping $\varphi: A\rightarrow E$
is Lipschitz continuous. Then,
\begin{enumerate}[(a)]
\item If $\mathcal{H}^{m-k}(\varphi^{-1}\{\bz\})>0$ for $\mathcal{H}^{k}$ almost all $\bz$ in $E$,
then $E$ is $(\mathcal{H}^{k},k)$ rectifiable.
\item If $E$ is $(\mathcal{H}^{k},k)$ rectifiable, then the set $\varphi^{-1}\{\bz\}$ is $(\mathcal{H}^{m-k},m-k)$ rectifiable
for $\mathcal{H}^{k}$ almost all $\bz$, and
\[\int_{A}p(\by)J_{k}^\varphi(\by)\mathcal{H}^{m}(d\by)=\int_{E}\int_{
\varphi^{-1}\{\bz\}}p(\by)\mathcal{H}^{m-k}(d\by)\mathcal{H}^{k}(d\bz),\]where
$p(\cdot)$ is an $\mathcal{H}^{m}$ integrable function.
\end{enumerate}
\end{lemma}

The area formula can be viewed as a generalization of the concept of
change of variables in integration. When $n=m$, it is the typically
change-of-variables formula in calculus. This extension to cases
with $n>m$ provides a powerful tool for computing integrals taken
over $m$-dimensional surfaces in the $n$-dimensional space, which is
transformed into a Lebesgue integration over a subset of
$\mathbb{R}^m$. In particular, the $m$-dimensional area of the image
of a continuously differentiable mapping $u$ from a domain $E\subset
\mathbb{R}^{m}$ into $\mathbb{R}^{n}$ is defined as the integral of
the Jacobian $J_{m}^u$ over $E$. The coarea formula can be viewed as
an extension of the Fubini's Theorem to Hausdorff measure, i.e., a
double integral can be computed using iterated integrals.

\subsection{Proof of Theorem \ref{thm: expecttransform}}

\noindent{\it Proof.} Note that $J_{m}^u(\bx)\neq 0$ for almost all
$\bx\in \Omega$ up to a null set. Then,
\begin{eqnarray*}
  \E\left[p(\bX)\right] &=& \E\left[p(\bX)\right] = \int_{\Omega}p(\bx)f(\bx)d\bx=\int_{\Omega}\frac{p(\bx)f(\bx)}{J_{m}^u(\bx)}J_{m}^u(\bx)d\bx\\
  &=& \int_{u\{\Omega\}}\frac{p(u^{-1}(\by))f(u^{-1}(\by))}{J_{m}(u^{-1}(\by))}\mathcal{H}^{m}(d\by),
\end{eqnarray*}
where the last equality follows from the area formula in Lemma
\ref{lem:areaformula}.

%Suppose there exists a Lipschitz continuous mapping
%\[\varphi: u\{D\}\mapsto \mathbb{R}^{q},\]such that the image set $\varphi\{u\{\Omega\}\}$ is
%$(\mathcal{H}^{k}, k)$ rectifiable, where $k\le\min\{q,m\}$.

Note that by definition, $u\{\Omega\}$ is $({\cal H}^m, m)$
rectifiable because $u$ is a continuously differentiable mapping. If
$\varphi\{u\{\Omega\}\}$ is $({\cal H}^k,k)$ rectifiable, applying
the coarea formula in Lemma \ref{lem:coareaformula} yields
\begin{eqnarray}
 \lefteqn{\E\left[p(\bX)\right]=\int_{u\{\Omega\}}\frac{p(u^{-1}(\by))f(u^{-1}(\by))}{J_{m}^u(u^{-1}(\by))}\mathcal{H}^{m}(d\by)}\nonumber\\
 &=&\int_{u\{\Omega\}}\frac{p(u^{-1}(\by))f(u^{-1}(\by))J_{k}^\varphi(\by)}{J_{m}^u(u^{-1}(\by))J_{k}^\varphi(\by)}\mathcal{H}^{m}(d\by)\nonumber\\
 &=&\int_{\varphi\{u\{\Omega\}\}}\int_{\varphi^{-1}\{\bz\}}\frac{p(u^{-1}(\by))f(u^{-1}(\by))}{J_{m}^u(u^{-1}(\by))J_{k}^\varphi(\by)}
 \mathcal{H}^{m-k}(d\by)\mathcal{H}^{k}(d\bz),\label{eqn:ProofTemp3}
\end{eqnarray}
where the last equality follows from the coarea formula in Lemma
\ref{lem:coareaformula}.

For $w(\bz)$ defined in Equation (\ref{eqn:Weight}), it can be seen
that
\begin{eqnarray}
\lefteqn{\E\left[p(\bX)\right]=\int_{\varphi\{u\{\Omega\}\}}w(\bz)\int_{\varphi^{-1}(\bz)}\frac{f(u^{-1}(\by))}{J_{m}^u(u^{-1}(\by))J_{k}^\varphi(\by)}\mathcal{H}^{m-k}(d\by)\mathcal{H}^{k}(d\bz)}\nonumber\\
 &=&\int_{\varphi\{u\{\Omega\}\}}\int_{\varphi^{-1}\{\bz\}}w(\bz)\frac{f(u^{-1}(\by))}{J_{m}^u(u^{-1}(\by))J_{k}^\varphi(\by)}\mathcal{H}^{m-k}(d\by)\mathcal{H}^{k}(d\bz)\nonumber\\
 &=&\int_{u\{\Omega\}}w(\varphi(\by))\frac{f(u^{-1}(\by))}{J_{m}^u(u^{-1}(\by))}\mathcal{H}^{m}(d\by),\label{eqn:ProofTemp1}
\end{eqnarray}
where the last equality follows from the coarea formula in Lemma
\ref{lem:coareaformula}.

Then by the area formula in Lemma \ref{lem:areaformula},
\begin{eqnarray}
  \lefteqn{\E\left[p(\bX)\right]=\int_{u(\Omega)}w(\varphi(\by))\cdot\frac{f(u^{-1}(\by))}{J_{m}^u(u^{-1}(\by))}\mathcal{H}^{m}(d\by)}\nonumber\\
  &=&\int_{\Omega}w\left(\varphi(u(\bx))\right){f(\bx)}d\bx=\E\left[w\left(\varphi(u(\bX))\right)\right].\label{eqn:ProofTemp2}
  \end{eqnarray}

The rest of proof will be devoted to the second half, i.e.,
$\Var\left[w\left(\varphi(u(\bX))\right)\right]\le
\Var\left[p(\bX)\right]$. Because these two integrands have the same
mean, it suffices to show that
\[\E\left[w^2\left(\varphi(u(\bX))\right)\right]\le
\E\left[p^2(\bX)\right].\]Note that for any given $\bz\in
\varphi\left\{u\{\Omega\}\right\}$,
\begin{eqnarray*}
\lefteqn{w^2(\bz)=\left.\left[\int_{\varphi^{-1}\{\bz\}}\frac{p(u^{-1}(\by))f(u^{-1}(\by))}{J_{m}^u(u^{-1}(\by))J_{k}^\varphi(\by)}\mathcal{H}^{m-k}(d\by)\right]^2\right/\left[\int_{\varphi^{-1}\{\bz\}}\frac{f(u^{-1}(\by))}{J_{m}^u(u^{-1}(\by))J_{k}^\varphi(\by)}\mathcal{H}^{m-k}(d\by)\right]^2}\\
&\le
&{\int_{\varphi^{-1}\{z\}}\frac{p^2(u^{-1}(\by))f(u^{-1}(\by))}{J_{m}^u(u^{-1}(\by))J_{k}^\varphi(\by)}\mathcal{H}^{m-k}(d\by)\int_{\varphi^{-1}\{\bz\}}\frac{f(u^{-1}(\by))}{J_{m}^u(u^{-1}(\by))J_{k}^\varphi(\by)}\mathcal{H}^{m-k}(d\by)\over\left[\int_{\varphi^{-1}\{\bz\}}\frac{f(u^{-1}(\by))}{J_{m}^u(u^{-1}(\by))J_{k}^\varphi(\by)}\mathcal{H}^{m-k}(d\by)\right]^2}\\
&=&
\left.\int_{\varphi^{-1}\{\bz\}}\frac{p^2(u^{-1}(\by))f(u^{-1}(\by))}{J_{m}^u(u^{-1}(\by))J_{k}^\varphi(\by)}\mathcal{H}^{m-k}(d\by)\right/\int_{\varphi^{-1}\{\bz\}}\frac{f(u^{-1}(\by))}{J_{m}^u(u^{-1}(\by))J_{k}^\varphi(\by)}\mathcal{H}^{m-k}(d\by),
\end{eqnarray*}
where the inequality follows from Cauchy-Schwarz inequality
$\left|\int q_1(\by)q_2(\by)\,d\by\right|^2\le \int
q_1^2(\by)\,d\by\int q_2^2(\by)\,d\by$ by letting
\[q_1(\by)=p(u^{-1}(\by))\sqrt{f(u^{-1}(\by))}/\sqrt{J_{m}^u(u^{-1}(\by))J_{k}^\varphi(\by)},\quad
q_2(\by)=\sqrt{f(u^{-1}(\by))}/\sqrt{J_{m}^u(u^{-1}(\by))J_{k}^\varphi(\by)}.\]

Then similar to the arguments in Equations (\ref{eqn:ProofTemp1})
and (\ref{eqn:ProofTemp2}), we have
\begin{eqnarray*}
\lefteqn{\E\left[w^2\left(\varphi(u(\bX))\right)\right]=\int_{\varphi\{u\{\Omega\}\}}w^2(\bz)\int_{\varphi^{-1}(\bz)}\frac{f(u^{-1}(\by))}{J_{m}^u(u^{-1}(\by))J_{k}^\varphi(\by)}\mathcal{H}^{m-k}(d\by)\mathcal{H}^{k}(d\bz)}\\
&\le
&\int_{\varphi\{u\{\Omega\}\}}\int_{\varphi^{-1}\{\bz\}}\frac{p^2(u^{-1}(\by))f(u^{-1}(\by))}{J_{m}^u(u^{-1}(\by))J_{k}^\varphi(\by)}\mathcal{H}^{m-k}(d\by)\mathcal{H}^{k}(d\bz)\\
&=& \E\left[p^2(\bX)\right],
\end{eqnarray*}
where the last equality follows from a similar argument as in
Equation (\ref{eqn:ProofTemp3}), and the proof is completed.
\hfill\blot

\subsection{Proof of Theorem \ref{thm:Sensitivity}}

{\noindent \it Proof.} Note that by assumption, ${\cal
H}^1(\varphi^{-1}\{\bz\})>0$ for almost all $\bz$. Then by the first
part of Lemma \ref{lem:coareaformula}, $\varphi\{u\{\Omega\}\}$ is
$({\cal H}^{m-1},m-1)$ rectifiable.

Applying Theorem \ref{thm: expecttransform}, we have
\[\alpha(\xi)=\E\left[w\left(\varphi(u(\bX));\xi\right)\right].\]

Because $w\left(\varphi(u(\bX));\xi\right)$ satisfies the Lipschitz
continuous condition w.r.t. $\xi$, by the dominated convergence
theorem (Durrett 2005), we interchange the order of differentiation
and expectation, leading to
\[\alpha'(\xi)=\E\left[\partial_{\xi}w\left(\varphi(u(\bX));\xi\right)\right]=\E\left[\nu\left(\varphi(u(\bX));\xi\right)\right],\]which
completes the proof. \hfill\blot

%\subsection{A Lemma of Liu and Hong (2011)}\label{sec:LemmaLiuHong}

%\begin{assumption}\label{assum:Lip}
%For any $\theta\in \Theta$, $l(\bX(\theta))$ and $h(\bX(\theta))$
%are differentiable with respect to $\theta$ with probability 1
%(w.p.1), and there exist random variables $K_l$ and $K_h$ with
%finite second moments that may depend on $\theta$, such that
%$|l(\bX(\theta+\Delta\theta))-l(\bX(\theta))|\le K_l|\Delta\theta|$
%and $|h(\bX(\theta+\Delta\theta))-h(\bX(\theta))|\le
%K_h|\Delta\theta|$ when $\Delta\theta$ is small enough.
%\end{assumption}

%\begin{assumption}\label{assum:cont}
%For any $\theta\in\Theta$, $\partial_\theta v(\theta,u)$ exists and
%is continuous at $(\theta,y)$, where
%$v(\theta,u)=\E\left[l(\bX)\cdot1_{\{h(\bX)\le u\}}\right]$.
%\end{assumption}

%Liu and Hong (2011) prove the following result.

%\begin{lemma}[Liu and Hong 2011]\label{thm:LiuHong11}
%Suppose that $\E\left[|l(\bX(\theta))|^2\right]<\infty$,
%$\E\left[|h(\bX(\theta))|^2\right]<\infty$, and Assumptions
%\ref{assum:Lip} and \ref{assum:cont} are satisfied. Then,
%\begin{equation*}
%\alpha'(\theta) = \E\left[\partial_\theta
%l(\bX(\theta))\cdot1_{\{h(\bX(\theta))\le \xi\}}\right] -
%\partial_{\xi} \E\left[l(\bX(\theta))\partial_\theta h(\bX(\theta))\cdot1_{\{h(\bX(\theta))\le
%\xi\}}\right].
%\end{equation*}
%\end{lemma}

\newpage

\section{Online Supplement}

\subsection{Derivation of Estimators for Section
\ref{GreekExperiments}}\label{sec:OnlineDerive}

\subsubsection{Derivation under the Black-Scholes
Model}\label{sec:OnlineBS} Under the Black-Scholes (BS) model,
\[X_{i+1}=X_i\exp\left((r-\sigma^2/2)T/m+\sigma \sqrt{T/m}N_{i+1}\right),\quad i=0,\dots,m-1,\]where
$\{N_1,\dots,N_m\}$ are independent standard normal random
variables, and the initial underlying asset price $X_0=x_0$ is a
constant.

To simply notation, we let $\bar X = \sum_{i=1}^{m}X_i/m$, $\widehat
X = \max(X_1,\dots,X_m)$, $\mu = r-\sigma^2/2$, and $\tau = T/m$.
For $i=1,\dots,m-1$, the conditional density of $X_{i+1}$ given
$X_i=x_i$ is
\begin{equation}\label{eqn:TransitionBS}
f_{i+1}(x_{i+1}|x_i)={1\over \sigma
x_{i+1}\sqrt{\tau}}\phi\left({1\over
\sigma\sqrt{\tau}}\left(\log(x_{i+1}/x_i)-\mu \tau\right)\right),
\end{equation}
where $\phi$ denotes the standard
normal density function. Let $\bx=(x_1,\dots,x_m)$. The joint
density of $\bX=(X_1,\ldots,X_m)$ is
\[f(\bx) = \prod_{i=0}^{m-1}f_{i+1}(x_{i+1}|x_i).\]

\noindent{\it \textbf{The Change-of-Variables Approach}}

We apply Proposition \ref{pro:main} to derive change-of-variables
estimators, where $h(\bx)$ shall be specified in the context.

\begin{itemize}
\item \textbf{The digital option}.

Note that $dX_m/dx_0=X_m/x_0$. By Theorem 1 of Liu and Hong (2011),
\begin{eqnarray*}
{\tt delta}={d\over dx_0}\E\left[e^{-rT}1_{\{X_m\ge
K\}}\right]=-{d\over dK}\E\left[e^{-rT}{dX_m\over dx_0}1_{\{X_m\ge
K\}}\right]={d\over dK}\E\left[e^{-rT}{X_m\over x_0}1_{\{X_m\le
K\}}\right].
\end{eqnarray*}
Then setting $h(\bx)=x_m$, $g(\bx)=e^{-rT}x_m/x_0$, and applying
Proposition \ref{pro:main}, we have $\bar J_m(\bz,t)=1/t^{m-1}$, and
\begin{equation*}
{\tt delta}=\E\left[\nu(\bZ;K)\right]=\E\left[e^{-rT}{X_1K\over
X_mx_0}f_1\left(X_1K/X_m|x_0\right)\right],
\end{equation*}
where $\bZ\triangleq (Z_1,\dots,Z_m)=\bX/h(\bX)$, and for
$\bz=(z_1,\dots,z_m)$,
\begin{eqnarray*}
\nu(\bz;\xi)&=&g(\bz \xi)f(\bz \xi)\xi^{m-1}\left/\int_{0}^\infty
f(t\bz)t^{m-1}dt\right.\\
&=& e^{-rT}z_m\xi/x_0f_1(z_1 \xi)\left/\int_{0}^{\infty}f_1(z_1 t)\,dt\right.\\
&=&{e^{-rT}z_mz_1\xi\over x_0}f_1(z_1\xi)={e^{-rT}z_1\xi\over
x_0}f_1(z_1\xi),
\end{eqnarray*}
where the last equality follows from $z_m=1$.

In a similar manner, it can be derived that
\begin{eqnarray*}
{\tt theta}&=& \E\left[re^{-rT}1_{\{X_m\ge
K\}}\right]-\E\left[e^{-rT}{X_1K\over X_m}f_1\left(X_1K/
X_m|x_0\right)(\mu T+\log(K/x_0))/(2T)\right],\\
{\tt vega}&=& \E\left[e^{-rT}{X_1K\over X_m}f_1\left(X_1K/X_m|x_0\right)(\log(K/x_0)-(\mu+\sigma^2)T)/\sigma\right],\\
{\tt gamma}&=& {d{\tt delta}\over dx_0}=\E\left[{e^{-rT}X_1K\over
x_0^2\sigma^2\tau X_m}f_1\left(X_1K/
X_m|x_0\right)\left(\log(X_1K/x_0X_m)-(\mu+\sigma^2)\tau\right)\right].
\end{eqnarray*}
where $f_1$ is the conditional density of $X_1$ given $x_0$ as
specified in (\ref{eqn:TransitionBS}).

\item \textbf{The Asian digital option}.

The derivation is parallel to that for the digital option, except
that $h(\bx)=\sum_{i=1}^{m}x_m/m$, and $\bar
J_m(\bz,t)=\sqrt{m}/t^{m-1}$. Applying Proposition \ref{pro:main},
we have
\begin{eqnarray*}
{\tt delta}&=& \E\left[{e^{-rT}X_1K \over \bar X}f_1\left(X_1K/\bar
X|x_0\right)\right],\\
{\tt theta}&=& \E\left[re^{-rT}1_{\{\bar X\ge
K\}}\right]-\E\left[{e^{-rT}X_1K\over 2mT\bar X}f_1\left(X_1K/\bar
X|x_0\right)\sum_{i=1}^{m}{X_i\over \bar X}\left(\log{X_iK\over
x_0\bar
X}-\mu_bi\tau\right)\right],\\
{\tt vega}&=& \E\left[{e^{-rT}X_1K\over m\sigma\bar
X}f_1\left(X_1K/\bar X|x_0\right)\sum_{i=1}^{m}{X_i\over \bar
X}\left(\log{X_iK\over x_0\bar
X}-(\mu+\sigma^2)i\tau\right)\right],\\
{\tt gamma}&=&\E\left[{e^{-rT}X_1K\over x_0^2\sigma^2\tau \bar
X}f_1\left(X_1K/\bar X|x_0\right)\left(\log{ X_1K\over x_0\bar
X}-(\mu+\sigma^2)\tau\right)\right].
\end{eqnarray*}

\item \textbf{The barrier call option}.

Note that $dX_i/dx_0=X_i/x_0$ for $i=1,\dots,m$. By Theorem 1 of Liu
and Hong (2011),
\begin{eqnarray*}
\lefteqn{{\tt delta}={d\over
dx_0}\E\left[e^{-rT}(X_m-K)^+1_{\{\widehat X\le
\kappa\}}\right]}\\
&=& \E\left[e^{-rT}{X_m\over x_0}1_{\{X_m\ge K\}}1_{\{\widehat X\le
\kappa\}}\right]-{d\over d\kappa}\E\left[e^{-rT}(X_m-K)^+{\widehat
X\over x_0}1_{\{\widehat X\le \kappa\}}\right].
\end{eqnarray*}
Similar to the derivation for the digital option, we apply
Proposition \ref{pro:main} to derive change-of-variables estimator
for the second term on the right-hand-side of the above equation. In
particular, we set $h(\bx)=\widehat \bx\triangleq
\max(x_1,\dots,x_m)$, and $i^* = \textrm{argmax} (X_1,\dots, X_m)$.
Then the Jacobian $\bar J_m(\bz,t)=1/t^{m-1}$, and Proposition
\ref{pro:main} leads to
\[{d\over d\kappa}\E\left[e^{-rT}(X_m-K)^+{\widehat
X\over x_0}1_{\{\widehat X\le
\kappa\}}\right]=\E\left[{e^{-rT}\kappa X_1\over x_0\widehat
X}f_1\left(X_1\kappa/\widehat X|x_0\right)\left(X_m\kappa/\widehat
X-K\right)^+\right],\]and thus
\begin{eqnarray*}
{\tt delta}&=& \E\left[e^{-rT}{X_m\over x_0}1_{\{X_m\ge
K\}}1_{\{\widehat X\le \kappa\}}-{e^{-rT}\kappa X_1\over x_0\widehat
X}f_1\left(X_1\kappa/\widehat X|x_0\right)\left(X_m\kappa/\widehat
X-K\right)^+\right].
\end{eqnarray*}
Similarly it can be derived that
\begin{eqnarray*}
{\tt theta}&=&\E\left[re^{-rT}(X_m-K)^+1_{\{\widehat X\le
\kappa\}}\right]-{e^{-rT}\over 2T} \E\left[X_m1_{\{X_m\ge K\}}1_{\{\widehat
X\le \kappa\}}(\log(X_m/x_0)+\mu T)\right]\\
&&+{e^{-rT}\over 2T}\E\left[{X_1\kappa\over \widehat
X}f_1\left(X_1\kappa/\widehat X|x_0\right)\left(X_m\kappa/\widehat
X-K\right)^+(\log(\kappa/x_0)+\mu i^*\tau)\right],\\
{\tt vega}&=& \E\left[e^{-rT}X_m1_{\{X_m\ge K\}}1_{\{\widehat X\le
\kappa\}}\left(\log(X_m/x_0)-(\mu+\sigma^2)T\right)/\sigma\right]\\
&&-\E\left[{e^{-rT}X_1\kappa\over \sigma\widehat
X}f_1\left(X_1\kappa/\widehat X|x_0\right)\left(X_m\kappa/\widehat
X-K\right)^+\left(\log(\kappa/x_0)-(\mu+\sigma^2)i^*\tau\right)\right],\\
{\tt gamma}&=& E\left[{e^{-rT}K^2 X_1 f_1(X_1K/X_m |x_0)\over X_m x_{0}^{2}} 1_{\{\widehat
XK/X_m\leq \kappa \}} \right]\\
&&-\E\left[{e^{-rT}X_1 X_m\kappa^2 \over x_0^2 \widehat
X^2}f_1(X_1 \kappa/\widehat X|x_0)1_{\{X_m\kappa/\widehat
X\geq K\}}\right]\\
&&-\E\left[{e^{-rT}X_1\kappa\over
\widehat X x_0^2\sigma^2\tau}f_1(X_1 \kappa /\widehat X)(X_m\kappa/\widehat
X-K)^+ \left(\log{X_1\kappa\over
x_0\widehat X}-(\mu+\sigma^2)\tau\right)\right].
\end{eqnarray*}

\end{itemize}

\noindent{\it \textbf{The Likelihood Ratio Approach}}

Note that the likelihood ratios with respect to $S_0$, $\sigma$ and
$T$ are
\begin{eqnarray*}
L_1(\bx)&=&{d\over dx_0}\log f(\bx)={1\over
x_0\sigma^2\tau} \left(\log {x_1\over x_0}-\mu\tau\right),\\
L_2(\bx)&=&{1\over f(\bx)}{d^2\over dx_0^2}
f(\bx)=L_1^2(\bx)-{1\over x_0^2\sigma^2\tau}\left(\log {x_1\over x_0}-\mu\tau+1\right),\\
L_3(\bx)&=&{d\over d T}\log f(\bx)={1\over m}\sum_{i=0}^{m-1}\left[{1\over 2\sigma^2\tau^2}\left(\log {x_{i+1}\over x_i}-\mu\tau\right)^2+{\mu\over \sigma^2\tau}\left(\log {x_{i+1}\over x_i}-\mu\tau\right)-{1\over 2\tau}\right],\\
L_4(\bx)&=&{d\over d\sigma}\log f(\bx)=\sum_{i=0}^{m-1}\left[{1\over
\sigma^3\tau}\left(\log {x_{i+1}\over
x_i}-\mu\tau\right)^2-{1\over \sigma}\left(\log {x_{i+1}\over
x_i}-\mu\tau\right)-{1\over \sigma}\right].
\end{eqnarray*}

\begin{itemize}
\item \textbf{The digital option}.
\begin{eqnarray*}
\text{\tt delta}&=&{d\over dx_0}\E\left[e^{-rT}1_{\{X_m\ge
K\}}\right]=\E\left[e^{-rT}L_1(\bX)1_{\{X_m\ge
K\}}\right],\\
\text{\tt theta}&=&-{d\over dx_0}\E\left[e^{-rT}1_{\{X_m\ge
K\}}\right]=\E\left[re^{-rT}1_{\{X_m\ge
K\}}\right]-\E\left[e^{-rT}L_3(\bX)1_{\{X_m\ge
K\}}\right],\\
\text{\tt vega}&=&{d\over d\sigma}\E\left[e^{-rT}1_{\{X_m\ge
K\}}\right]=\E\left[e^{-rT}L_4(\bX)1_{\{X_m\ge
K\}}\right],\\
\text{\tt gamma}&=&{d^2\over dx_0^2}\E\left[e^{-rT}1_{\{X_m\ge
K\}}\right]=\E\left[e^{-rT}L_2(\bX)1_{\{X_m\ge K\}}\right].
\end{eqnarray*}
\item \textbf{The Asian digital option}.
\begin{eqnarray*}
\text{\tt delta}&=&{d\over dx_0}\E\left[e^{-rT}1_{\{\bar X\ge
K\}}\right]=\E\left[e^{-rT}L_1(\bX)1_{\{\bar X\ge
K\}}\right],\\
\text{\tt theta}&=&-{d\over dx_0}\E\left[e^{-rT}1_{\{\bar X\ge
K\}}\right]=\E\left[re^{-rT}1_{\{\bar X\ge
K\}}\right]-\E\left[e^{-rT}L_3(\bX)1_{\{\bar X\ge
K\}}\right],\\
\text{\tt vega}&=&{d\over d\sigma}\E\left[e^{-rT}1_{\{\bar X\ge
K\}}\right]=\E\left[e^{-rT}L_4(\bX)1_{\{\bar X\ge
K\}}\right],\\
\text{\tt gamma}&=&{d^2\over dx_0^2}\E\left[e^{-rT}1_{\{\bar X\ge
K\}}\right]=\E\left[e^{-rT}L_2(\bX)1_{\{\bar X\ge K\}}\right].
\end{eqnarray*}
\item \textbf{The barrier call option}.
\begin{eqnarray*}
\text{\tt delta}&=&{d\over
dx_0}\E\left[e^{-rT}(X_m-K)^+1_{\{\widehat X\le
\kappa\}}\right]=\E\left[e^{-rT}L_1(\bX)(X_m-K)^+1_{\{\widehat X \le
\kappa\}}\right],\\
\text{\tt theta}&=&-{d\over
dx_0}\E\left[e^{-rT}(X_m-K)^+1_{\{\widehat X \le
\kappa\}}\right]\\
&=&\E\left[re^{-rT}(X_m-K)^+1_{\{\widehat X\le
\kappa\}}\right]-\E\left[e^{-rT}L_3(\bX)(X_m-K)^+1_{\{\widehat X\le
\kappa\}}\right],\\
\text{\tt vega}&=&{d\over
d\sigma}\E\left[e^{-rT}(X_m-K)^+1_{\{\widehat X\le
\kappa\}}\right]=\E\left[e^{-rT}L_4(\bX)(X_m-K)^+1_{\{\widehat X\le
\kappa\}}\right],\\
\text{\tt gamma}&=&{d^2\over
dx_0^2}\E\left[e^{-rT}(X_m-K)^+1_{\{\widehat X\le
\kappa\}}\right]=\E\left[e^{-rT}L_2(\bX)(X_m-K)^+1_{\{\widehat X\le
\kappa\}}\right].
\end{eqnarray*}
\end{itemize}

\noindent{\it \textbf{Conventional Conditional Monte Carlo
Approach}}
\begin{itemize}
\item \textbf{The digital option}.

Conditioning on $X_{m-1}$ yields
\[\E\left[e^{-rT}1_{\{X_m\ge
K\}}\right]=\E\left(\E\left[e^{-rT}1_{\{X_m\ge
K\}}|X_{m-1}\right]\right)=\E\left[e^{-rT}\left(1-\Phi\left({1\over\sigma\sqrt{\tau}}\left(\log{K\over
X_{m-1}}-\mu\tau\right)\right)\right)\right],\]where $\Phi$
denotes the standard normal distribution function.

Using the pathwise method on the right-hand-side (RHS) of the above
equation, we have
\begin{eqnarray*}
{\tt delta}&=&\E\left[{e^{-rT}K f_m(K|X_{m-1})\over
x_0}\right],\\
{\tt theta}&=&\E\left[re^{-rT}1_{\{X_m\ge
K\}}\right]-{1\over 2T}\E\left[{e^{-rT}K f_m(K|X_{m-1})}\left({\log(K/x_0)+\mu T}\right)\right],\\
{\tt vega}&=&\E\left[{e^{-rT}K f_m(K|X_{m-1})\left(\log(K/x_0)-(\mu+\sigma^2)T\right)\over\sigma}\right],\\
{\tt gamma}&=& \E\left[{e^{-rT}Kf_m(K|X_{m-1})(\log(K/X_{m-1})-(\mu+\sigma^2)\tau)\over
x_0^2\sigma^2\tau}\right].
\end{eqnarray*}

\item \textbf{The Asian digital option}.

Let $S_{-m}=\sum_{i=1}^{m-1}X_i$ and $S_{m}=mK-S_{-m}$. Conditioning on
$\{X_1,\dots,X_{m-1}\}$, we have
\begin{eqnarray*}
\E\left[e^{-rT}1_{\{\bar X\ge
K\}}\right]&=&\E\left[\E\left[e^{-rT}1_{\{\bar X\ge
K\}}|X_1,\dots,X_{m-1}\right]\right]\\
&=&\E\left[e^{-rT}\left(1-\Phi\left({1\over\sigma\sqrt{\tau}}\left(\log{mK-S_{-m}\over
X_{m-1}}-\mu\tau\right)\right)\right)1_{\{S_{-m}\le mK\}}\right].
\end{eqnarray*}
Using the pathwise method on the RHS of the above equation, we have
\begin{eqnarray*}
{\tt delta}&=&\E\left[{e^{-rT}mK f_m(S_m|X_{m-1})\over
x_0}1_{\{S_{-m}< mK\}}\right],\\
{\tt theta}&=&\E\left[re^{-rT}1_{\{\bar X\ge
K\}}\right]-{1\over 2T}\E\left[{e^{-rT}}f_m(S_m|X_{m-1})1_{\{S_{-m}< mK\}}\right.\\
&&\left.{S_m(\log (S_{m}/x_0)+\mu T)}+ \sum_{i=1}^{m-1}{X_i(\log(X_i/x_0)+i\mu\tau)}\right],\\
{\tt vega}&=&\E\left[{e^{-rT}f_m(S_m|X_{m-1})\over \sigma}1_{\{S_{-m}<
mK\}}\right.\\
&&\left.\left({S_m(\log({S_{m}/x_0})-(\mu+\sigma^2)T)}+\sum_{i=1}^{m-1}{X_i(\log(X_i/x_0)-i(\mu+\sigma^2)\tau)}\right)\right],\\
{\tt gamma}&=&\E\left[{e^{-rT}mKf_m(S_m|X_{m-1})\over S_m x_0^2}1_{\{S_{-m}<
mK\}}\left({S_{-m}- S_{m}}+{mK(\log{(S_{m}/
X_{m-1})}-\mu\tau)\over\sigma^2\tau}\right)\right].
\end{eqnarray*}

\end{itemize}

\subsubsection{Derivation under the Variance Gamma
Model}\label{sec:VG}

Under the variance gamma (VG) model,
\begin{equation*}
X_{i+1}=X_i\exp(\mu_{\beta} \tau +\theta G_{i+1}+\sigma
\sqrt{G_{i+1}}N_{i+1}), \quad i = 0,\dots,m-1,
\end{equation*}
where $\tau=T/m$, $\{G_1,\dots,G_m\}$ are independent gamma random
variables with scale parameter $\tau/\beta$ and shape parameter
$\beta$, $\{N_1,\dots,N_m\}$ are independent standard normal random
variables, and $\mu_{\beta}=r + 1/\beta\log
(1-\theta\beta-\sigma^2\beta/2)$.

Note that for $i=0,\dots,m-1$, the conditional density of $X_{i+1}$
given $G_{i+1}$ and $X_i=x_i$ is
\[f_{i+1}(x_{i+1}|x_i,G_{i+1})={1\over x_{i+1}\sigma\sqrt{G_{i+1}}}\phi\left({1\over \sigma\sqrt{G_{i+1}}}\left(\log\left({x_{i+1}\over x_i}\right)-\mu_{\beta} \tau-\theta
G_{i+1}\right)\right),\]and thus the conditional density of
$\bX=(X_1,\dots,X_m)$ given $\bG=(G_1,\dots,G_m)$ is
\[f(\bx|\bG)=\prod_{i=0}^{m-1}f_{i+1}(x_{i+1}|x_i,G_{i+1})=\prod_{i=0}^{m-1}{1\over x_{i+1}\sigma\sqrt{G_{i+1}}}\phi\left({1\over \sigma\sqrt{G_{i+1}}}\left(\log\left({x_{i+1}\over x_i}\right)-\mu_{\beta} \tau-\theta
G_{i+1}\right)\right),\]where $\phi$ denotes the standard normal
density function.

\noindent{\it \textbf{The Change-of-Variables Approach}}

The derivation of the change-of-variables estimators is similar to
that under the Black-Scholes model, except that we are working with
the conditional density of $\bX$ given $\bG$, instead of the
unconditional density. Specifically, for a parameter $\eta$ that
$\bX$ may depend on and a function $l(\cdot)$, it holds that
\begin{equation}\label{eqn:OnlineVG}
{d\over d\eta}\E\left[l(\bX)\right]={d\over
d\eta}\E\left(\E\left[l(\bX)|\bG\right]\right)=\E\left({d\over
d\eta}\E\left[l(\bX)|\bG\right]\right),
\end{equation}
provided that $\E\left[l(\bX)|\bG\right]$ is a smooth function of
$\bG$. For estimating price sensitivities under the VG model, it can
be verified that $\E\left[l(\bX)|\bG\right]$ is indeed smooth for
the function $l$ of concern, and therefore (\ref{eqn:OnlineVG}) is
justified. Then, Proposition \ref{pro:main} can be applied to
$d\E\left[l(\bX)|\bG\right]/d\eta$, leading to
\begin{equation}\label{eqn:OnlineVG1}
{d\over
d\eta}\E\left[l(\bX)|\bG\right]=\E\left[\psi(\bX,\bG)|\bG\right],
\end{equation}
for an appropriate function $\psi$, and thus combining with
(\ref{eqn:OnlineVG}) yields
\[{d\over d\eta}\E\left[l(\bX)\right]=\E\left(\E\left[\psi(\bX,\bG)|\bG\right]\right)=\E\left[\psi(\bX,\bG)\right].\]

Note that conditional on $\bG$, the VG model has the same structure
as the BS model. Therefore, under the VG model, the way of applying
Proposition \ref{pro:main} to derive $\psi$ in (\ref{eqn:OnlineVG1})
for various options is exactly parallel to that under the BS model
in Section \ref{sec:OnlineBS}. The details are thus omitted.

To further simplify notations, we let $S_G = \sum_{i=1}^{m}G_i$,
$\gamma=1-\theta\beta-\sigma^2\beta/2$, and
\[\phi_{1,d}=f_1\left(\left.{X_1K\over
X_m}\right|x_0,G_1\right),\quad \phi_{1,a}=f_1\left(\left.{X_1K\over
\bar X}\right|x_0,G_1\right),\quad
\phi_{1,b}=f_1\left(\left.{X_1\kappa\over \widehat
X}\right|x_0,G_1\right).\]

\begin{itemize}
\item \textbf{The digital option}.
\begin{eqnarray*}
{\tt delta}&=& \E\left[{e^{-rT}KX_1\over
x_0X_m}\phi_{1,d}\right],\\
{\tt theta}&=&\E\left[re^{-rT}1_{\{X_m\ge
K\}}\right]-\E\left[{e^{-rT}KX_1\phi_{1,d}\over 2TX_m}\left(\theta
S_G+\mu_{\beta}
T + \log{K\over x_0}\right)\right],\\
{\tt vega}&=& \E\left[{e^{-rT}KX_1\phi_{1,d}\over
X_m}\left({\log(K/x_0)-\mu_{\beta} T-\theta S_G\over
\sigma}-{T\sigma\over
\gamma}\right)\right],\\
{\tt gamma}&=& \E\left[{e^{-rT}KX_1\phi_{1,d}\over X_m x_0^2\sigma^2
G_1}\left(\log{KX_1\over x_0 X_m}-(\mu_{\beta}\tau+\theta
G_1+\sigma^2G_1)\right)\right].
\end{eqnarray*}
\item \textbf{The Asian digital option}.
\begin{eqnarray*}
{\tt delta}&=& \E\left[{e^{-rT}KX_1\over x_0 \bar
X}\phi_{1,a}\right],\\
{\tt theta}&=& \E\left[re^{-rT}1_{\{\bar X\ge
K\}}\right]-\E\left[{e^{-rT}KX_1\phi_{1,a}\over 2mT\bar
X}\sum_{i=1}^{m}\left(\log{KX_i\over x_0\bar X}+i\mu_{\beta}\tau+{X_i\over \bar X}\sum_{k=1}^{i}G_k \theta \right)\right],\\
{\tt vega}&=& \E\left[{e^{-rT}KX_1\over m\bar
X}\phi_{1,a}\sum_{i=1}^{m}{X_i\over \bar X}\left({1\over
\sigma}\left(\log{KX_i\over x_0\bar X}-i\mu_{\beta}\tau-\theta
S_G\right)-{i\sigma\tau\over \gamma}\right)\right],\\
{\tt gamma}&=& \E\left[{e^{-rT}KX_1\over \bar
X x_0^2\sigma^2G_1}\phi_{1,a}\left(\log{KX_1\over x_0 \bar X}-(\mu_{\beta}\tau+\theta
G_1+\sigma^2 G_1)\right)\right].
\end{eqnarray*}
\item \textbf{The barrier call option}.
Let $i^* = \textrm{argmax}\{X_1,\dots, X_m\}$.
\begin{eqnarray*}
{\tt delta}&=&\E\left[{e^{-rT}X_m\over x_0}1_{\{X_m\ge
K\}}1_{\{\widehat X\le \kappa\}}\right]- \E\left[{e^{-rT}\kappa
X_1\phi_{1,m}\over x_0\widehat X}(X_m\kappa/\widehat X-K)^+\right],\\
{\tt theta}&=&\E\left[re^{-rT}(X_m-K)^+1_{\{\widehat X\le \kappa\}}\right] -\E\left[{e^{-rT}X_m\over 2T}\left(\log{X_m\over x_0}+\mu_{\beta} T+\theta S_G\right)1_{\{X_m\ge K\}}1_{\{\widehat X\le \kappa\}}\right]\\
&&+\E\left[{e^{-rT}\kappa X_1\phi_{1,b}\over 2T\widehat X}(X_m\kappa/\widehat X-K)^+\left(\log{\kappa\over x_0}+\mu_{\beta} i^*\tau+\theta \sum_{j=1}^{i^*} G_j\right)\right],\\
{\tt vega}&=& \E\left[e^{-rT}X_m\left({\log(X_m/x_0)-\mu_{\beta} T-\theta
S_G\over \sigma}-{\sigma T\over \gamma}\right)1_{\{X_m\ge
K\}}1_{\{\widehat X\le \kappa\}}\right]\\
&&-\E\left[{e^{-rT}\kappa X_1\phi_{1,m}\over \widehat
X}(X_m\kappa/\widehat X-K)^+\left({\log(\kappa/x_0)-\mu_{\beta}
i^*\tau-\theta \sum_{j=1}^{i^*}G_j\over\sigma}-{\sigma
i^*\tau\over\gamma}\right)\right],\\
{\tt gamma}&=&\E\left[{e^{-rT}K^2X_1\phi_{1,d}\over x_0^2
X_m}1_{\{K\widehat X/X_m\le
\kappa\}}\right]-\E\left[{e^{-rT}\kappa^2 X_1X_m\phi_{1,b}\over x_0^2\widehat
X^2}1_{\{\kappa X_m/\widehat X\ge K\}}\right]\\
&&-\E\left[{e^{-rT}\kappa X_1\phi_{1,b}\over \widehat X x_0^2\sigma^2
G_1}(X_m\kappa/\widehat X-K)^+\left(\log{\kappa X_1\over
x_0\widehat X}-(\mu_{\beta}\tau+\theta G_1+\sigma^2 G_1)\right)\right].
\end{eqnarray*}
\end{itemize}

\noindent{\it \textbf{The Likelihood Ratio Approach}}

For any function $l$, and a market parameter $\eta$, note that
\begin{eqnarray*}
\lefteqn{{d\over d\eta}\E\left[l(\bX)\right]={d\over
d\eta}\E\left(\E\left[l(\bX)|\bG\right]\right)}\\
&=&{d\over d\eta}\E\left[\int
l(\bx)f(\bx|\bG)\,d\bx\right]=\E\left(\E\left[l(\bX){d\log
f(\bX|\bG)\over d\eta}|\bG\right]\right)=\E\left[l(\bX){d\log
f(\bX|\bG)\over d\eta}\right],
\end{eqnarray*}
where the interchange of expectation and differentiation in the
third equality is usually valid, because the integration $\int
l(\bx)f(\bx|\bG)\,d\bx$ is usually continuous in $\eta$ even when
$l$ is discontinuous.

By elementary algebra, the conditional likelihood ratios with
respect to $S_0$, $\sigma$ and $T$ are
\begin{eqnarray*}
L_1(\bx|\bG)&=&{d\over dx_0}\log f(\bx|\bG)={1\over
x_0\sigma^2G_0} \left(\log {x_1\over x_0}-(\mu_{\beta}\tau+\theta G_1)\right),\\
L_2(\bx|\bG)&=&{1\over f(\bx|\bG)}{d^2\over dx_0^2}
f(\bx|\bG)=L_1^2(\bx|\bG)-{1\over x_0^2\sigma^2G_1}\left(\log {x_1\over x_0}-(\mu_{\beta}\tau + \theta G_1)+1\right),\\
L_3(\bx|\bG)&=&{d\over d T}\log f(\bx|\bG)\\
&=&-{m\over 2T}+\sum_{i=1}^{m}{\left(\log{x_i\over
x_{i-1}}-(\mu_{\beta}\tau+\theta G_i)\right)^2\over 2\sigma^2T
G_i}+{\left(\log{x_i\over x_{i-1}}-(\mu_{\beta}\tau+\theta
G_i)\right)\left({\mu_{\beta}\over m}+{\theta G_i\over
T}\right)\over x_0^2\sigma^2
G_1},\\
L_4(\bx|\bG)&=&{d\over d\sigma}\log f(\bx|\bG)=-{m\over
\sigma}+\sum_{i=1}^{m}{\left(\log{x_i\over x_{i-1}}-(\mu_{\beta}\tau+\theta
G_i)\right)^2\over \sigma^3 G_i}-{\left(\log{x_i\over
x_{i-1}}-(\mu_{\beta}\tau+\theta G_i)\right)\tau\over \sigma
G_i(1-\theta\beta-\sigma^2\beta/2)}.
\end{eqnarray*}

\begin{itemize}
\item \textbf{The digital option}.
\begin{eqnarray*}
\text{\tt delta}&=&\E\left[e^{-rT}L_1(\bX|\bG)1_{\{X_m\ge
K\}}\right],\\
\text{\tt theta}&=&\E\left[re^{-rT}1_{\{X_m\ge
K\}}\right]-\E\left[e^{-rT}L_3(\bX|\bG)1_{\{X_m\ge
K\}}\right],\\
\text{\tt vega}&=&\E\left[e^{-rT}L_4(\bX|\bG)1_{\{X_m\ge
K\}}\right],\\
\text{\tt gamma}&=&\E\left[e^{-rT}L_2(\bX|\bG)1_{\{X_m\ge
K\}}\right].
\end{eqnarray*}
\item \textbf{The Asian digital option}.
\begin{eqnarray*}
\text{\tt delta}&=&\E\left[e^{-rT}L_1(\bX|\bG)1_{\{\bar X\ge
K\}}\right],\\
\text{\tt theta}&=&\E\left[re^{-rT}1_{\{\bar X\ge
K\}}\right]-\E\left[e^{-rT}L_3(\bX|\bG)1_{\{\bar X\ge
K\}}\right],\\
\text{\tt vega}&=&\E\left[e^{-rT}L_4(\bX|\bG)1_{\{\bar X\ge
K\}}\right],\\
\text{\tt gamma}&=&\E\left[e^{-rT}L_2(\bX|\bG)1_{\{\bar X\ge
K\}}\right].
\end{eqnarray*}
\item \textbf{The barrier call option}.
\begin{eqnarray*}
\text{\tt delta}&=&\E\left[e^{-rT}L_1(\bX|\bG)(X_m-K)^+1_{\{\widehat
X \le
\kappa\}}\right],\\
\text{\tt theta}&=&\E\left[re^{-rT}(X_m-K)^+1_{\{\widehat X\le
\kappa\}}\right]-\E\left[e^{-rT}L_3(\bX|\bG)(X_m-K)^+1_{\{\widehat
X\le
\kappa\}}\right],\\
\text{\tt vega}&=&\E\left[e^{-rT}L_4(\bX|\bG)(X_m-K)^+1_{\{\widehat
X\le
\kappa\}}\right],\\
\text{\tt gamma}&=&\E\left[e^{-rT}L_2(\bX|\bG)(X_m-K)^+1_{\{\widehat
X\le \kappa\}}\right].
\end{eqnarray*}
\end{itemize}

\noindent{\it \textbf{Conventional Conditional Monte Carlo
Approach}}

Define $\phi_{m}=f_m\left(\left.{K}\right|X_{m-1},G_m\right)$ in the following context.
\begin{itemize}
\item \textbf{The digital option}.

Conditioning on $\bG$ and $X_{m-1}$ and then applying the pathwise
method lead to
\begin{eqnarray*}
{\tt delta}&=& \E\left[{e^{-rT}K \phi_m\over x_0}\right],\\
{\tt theta}&=& \E\left[re^{-rT}1_{\{X_m\ge
K\}}\right]\\
&&-\E\left[{e^{-rT}\left(\theta\sum_{i=1}^m G_i+\mu_{\beta}
T+\log(K/x_0)\right)K \phi_m \over 2T }\right],\\
{\tt vega}&=& \E\left[{e^{-rT}K \phi_m}\left({\log(K/x_0)-\mu_{\beta}
T-\theta\sum_{i=1}^{m}G_i\over\sigma}-{T\sigma\over
1-\theta\beta-\sigma^2\beta/2}\right)\right]\\
{\tt gamma}&=&\E\left[{e^{-rT}K \phi_m \left(\log{K\over
X_{m-1}}-(\mu_{\beta}\tau+\theta G_m+\sigma^2 G_m)\right)\over x_0^2
\sigma^2{G_{m}}}\right].
\end{eqnarray*}

\item \textbf{The Asian digital option}.

Let $S_{-m}=\sum_{i=1}^{m-1}X_i$, $S_m = mK-S_{-m}$, and $S_G=\sum_{i=1}^mG_i$. To
further simplify notation, we let
$\gamma=1-\theta\beta-\sigma^2\beta/2$,
and \[\widehat{\phi}_m = f_m(S_{m}|X_{m-1},G_m).\]
Conditioning on $\bG$ and $(X_1,\dots,X_{m-1})$ and then applying
the pathwise method lead to
\begin{eqnarray*}
{\tt delta}&=& \E\left[{e^{-rT}mK\widehat{\phi}_m\over
x_0}1_{\{S_{-m}< mK\}}\right],\\
{\tt theta}&=& \E\left[re^{-rT}1_{\{\bar X\ge
K\}}\right]-\E\left[{e^{-rT}}\widehat \phi_m 1_{\{S_{-m}<
mK\}}\right.\\
&&\left.\left(S_m{\log(S_m/x_0)+\mu_{\beta} T+\theta S_G\over
2T}+{\sum_{i=1}^{m-1}X_i\left(\log{X_i\over x_0}+i\mu_{\beta}\tau+\theta
\sum_{k=1}^{i}G_k\right)}\right)\right],\\
{\tt vega}&=& \E\left[{e^{-rT}}\widehat \phi_m 1_{\{S_{-m}<
mK\}}\right.\\
&&\left.\left(S_m\left({\log{S_m\over x_0}-\mu_{\beta} T-\theta S_G\over\sigma}-{T\sigma\over \gamma}\right)+{\sum_{i=1}^{m-1}X_i\left({\log{X_i\over x_0}-i\mu_{\beta}\tau-\theta \sum_{j=1}^{i} G_j \over \sigma}-{i\sigma\tau\over \gamma}\right)}\right)\right],\\
{\tt gamma}&=& \E\left[{e^{-rT}mK \widehat \phi_m\over x_0^2
}1_{\{S_{-m}<
mK\}}\left(-1 + {S_{-m}\over S_m} +{mK\left(\log{S_{m}\over
X_{m-1}}-(\mu_{\beta}\tau + \theta G_m)\right)\over\sigma^2 G_m S_m}\right)\right].
\end{eqnarray*}

\end{itemize}

\subsection{Gradient Estimators for the Chance Constrained Program
Example}\label{sec:OnlineCCP}

\subsubsection{Multivariate Normal Distribution}

Suppose that $\bX=(X_1,\ldots,X_m)$ follows a multivariate normal
distribution with mean zero and a covariance matrix $\Sigma$. Its
density function is
\[f(\bx)={1\over (2\pi)^{m/2}}|\Sigma|^{-1/2}e^{-{1\over 2}\bx^T\Sigma^{-1}\bx}.\]

We first derive a conventional CMC estimator. Let $\Sigma_1$ denote
the covariance matrix of $(X_1,t_2X_2+\ldots+t_mX_m)$, and set
$A\triangleq(a_{ij})_{1\le i\le 2,1\le j\le2}=\Sigma_1^{-1}$.
Conditioning on $(X_2,\ldots,X_m)$, we have
\begin{eqnarray*}
\lefteqn{\Pr\left(\bt^T\bX\le b\right)}\\
&=&
\E\left[\Phi\left(\sqrt{a_{11}}\left({b-\sum_{i=2}^{m}t_iX_i\over
t_1}+{a_{12}\over a_{11}}\sum_{k=2}^{m}t_kX_k\right)\right)\right],
\end{eqnarray*}
where $\Phi$ denotes the standard normal distribution function.
Then,
\begin{eqnarray*}
\lefteqn{\partial_{t_1}\Pr\left(\bt^T\bX\le
b\right)}\\
&=&-\E\left[{\sqrt{a_{11}}(b-\sum_{i=2}^{m}t_iX_i)\over
t_1^2}\phi\left(\sqrt{a_{11}}\left({b-\sum_{i=2}^{m}t_iX_i\over
t_1}+{a_{12}\over a_{11}}\sum_{k=2}^{m}t_kX_k\right)\right)\right],
\end{eqnarray*}
where $\phi$ denotes the standard normal density function.

Next we consider the change-of-variables estimator by
(\ref{eqn:EstCCP}). Define $\bZ=\bX/(\bt^T\bX)$. By
(\ref{eqn:EstCCP}),
\begin{eqnarray*}
\lefteqn{\partial_{t_1}\Pr\left(\bt^T\bX\le b\right)}\\
&=&-\E\left[\left.b|b|^{m-1}Z_1e^{-{1\over 2}\bZ^T\Sigma^{-1}\bZ
b^2}\right/\int_{-\infty}^{\infty}|y|^{m-1}e^{-{1\over
2}\bZ^T\Sigma^{-1}\bZ y^2}\,dy\right].
\end{eqnarray*}
Note that
\begin{eqnarray*}
\lefteqn{\int_{-\infty}^{\infty}|y|^{m-1}e^{-{1\over
2}\bZ^T\Sigma^{-1}\bZ y^2}\,dy=2\int_{0}^{\infty}y^{m-1}e^{-{1\over
2}\bZ^T\Sigma^{-1}\bZ y^2}\,dy}\\
&=& (2\pi)^{-m/2}|\Sigma|^{-1/2}\left({2\over
\bZ^T\Sigma^{-1}\bZ}\right)^{m/2}\int_{0}^{\infty}e^{-u}u^{m/2-1}\,du\\
&=&(2\pi)^{-m/2}|\Sigma|^{-1/2}\left({2\over
\bZ^T\Sigma^{-1}\bZ}\right)^{m/2}\Gamma(m/2),
\end{eqnarray*}
where the second to last equality follows from a change of variables
$u = \bZ^T\Sigma^{-1}\bZ y^2/2$, and the gamma function is defined
by $\Gamma(t)=\int_{0}^{\infty}y^{t-1}e^{-y}\,dy$.

Therefore,
\begin{eqnarray*}
\partial_{t_1}\Pr\left(\bt^T\bX\le b\right)=-{b|b|^{m-1}\over
\Gamma(m/2)2^{m/2}}\E\left[Z_1\left(\bZ^T\Sigma^{-1}\bZ\right)^{m/2}e^{-{1\over
2}\bZ^T\Sigma^{-1}\bZ b^2}\right]
\end{eqnarray*}

\subsubsection{Multivariate Student's t-Distribution}

Suppose that $\bX=(X_1,\ldots,X_m)$ follows a multivariate Student's
t-distribution with mean zero, a covariance matrix $\Sigma$, and $v$
degrees of freedom. Its density function is
\[f(\bx)={C\over\left(1+\bx^T\Sigma^{-1}\bx/v\right)^{(v+m)/2}},\]where
\[C={\Gamma((v+m)/2)\over
\Gamma(v/2)v^{m/2}\pi^{m/2}|\Sigma|^{1/2}}.\]Equivalently, one may
generate a multivariate normal distribution $\bY=(Y_1,\ldots,Y_m)$
with mean zero and covariance $\Sigma$ and a chi square random
variable $\chi_v^2$ with $v$ degrees of freedoms, and set
$\bX=\bY\sqrt{v/\chi_v^2}$.

We first derive a conventional CMC estimator. Let $\Sigma_1$ denote
the covariance matrix of $(Y_1,t_2Y_2+\ldots+t_mY_m)$, and set
$A\triangleq(a_{ij})_{1\le i\le 2,1\le j\le2}=\Sigma_1^{-1}$.
Conditioning on $(Y_2,\ldots,Y_m)$ and $\chi_v^2$, we have
\begin{eqnarray*}
\lefteqn{\Pr\left(\bt^T\bX\le b\right)=\Pr\left(\bt^T\bY\le
b/\sqrt{v/\chi_v^2}\right)}\\
&=&
\E\left[\Phi\left(\sqrt{a_{11}}\left({b/\sqrt{v/\chi_v^2}-\sum_{i=2}^{m}t_iY_i\over
t_1}+{a_{12}\over a_{11}}\sum_{k=2}^{m}t_kY_k\right)\right)\right].
\end{eqnarray*}
Then,
\begin{eqnarray*}
\lefteqn{\partial_{t_1}\Pr\left(\bt^T\bX\le
b\right)}\\
&=&-\E\left[{\sqrt{a_{11}}(b/\sqrt{v/\chi_v^2}-\sum_{i=2}^{m}t_iY_i)\over
t_1^2}\phi\left(\sqrt{a_{11}}\left({b/\sqrt{v/\chi_v^2}-\sum_{i=2}^{m}t_iY_i\over
t_1}+{a_{12}\over a_{11}}\sum_{k=2}^{m}t_kY_k\right)\right)\right].
\end{eqnarray*}

Next we consider the new CMC estimator by (\ref{eqn:EstCCP}). Define
$\bZ=\bX/(t^T\bX)$. By (\ref{eqn:EstCCP}),
\begin{eqnarray*}
\lefteqn{\partial_{t_1}\Pr\left(t^T\bX\le b\right)}\\
&=&-\E\left[\left.{b|b|^{m-1}Z_1\over\left(1+\bZ^T\Sigma^{-1}\bZ
b^2/v\right)^{(v+m)/2}}\right/
\int_{-\infty}^{\infty}{|y|^{m-1}\over \left(1+\bZ^T\Sigma^{-1}\bZ
y^2/v\right)^{(v+m)/2}}\,dy\right].
\end{eqnarray*}
Note that
\begin{eqnarray*}
\lefteqn{\int_{-\infty}^{\infty}{|y|^{m-1}\over
\left(1+\bZ^T\Sigma^{-1}\bZ
y^2/v\right)^{(v+m)/2}}\,dy=2\int_{0}^{\infty}{y^{m-1}\over
\left(1+\bZ^T\Sigma^{-1}\bZ y^2/v\right)^{(v+m)/2}}\,dy}\\
&=& {v^{m/2}\over
(v+m-1)^{m/2}\left(\bZ^T\Sigma^{-1}\bZ\right)^{m/2}}2\int_{0}^{\infty}{w^{m-1}\over
\left(1+w^2/(v+m-1)\right)^{(v+m)/2}}\,dw\\
&=&{v^{m/2}\over
(v+m-1)^{m/2}\left(\bZ^T\Sigma^{-1}\bZ\right)^{m/2}}{C_1\Gamma(m/2)\Gamma(v/2)(v+m-1)^{(m-1)/2}\over
\sqrt{\pi}\Gamma((v+m-1)/2)}\\
&=& {v^{m/2}\Gamma(m/2)\Gamma(v/2)\over
\Gamma((v+m)/2)\left(\bZ^T\Sigma^{-1}\bZ\right)^{m/2}},
\end{eqnarray*}
where
$C_1=\Gamma((v+m)/2)\left/\left[\sqrt{(v+m-1)\pi}\Gamma((v+m-1)/2)\right]\right.$,
and the second to last equality follows from the closed form formula
of the raw moment for t-distribution.

Therefore,
\begin{eqnarray*}
\partial_{t_1}\Pr\left(\bt^T\bX\le b\right)=-\E\left[{b|b|^{m-1}Z_1\Gamma((v+m)/2)\left(\bZ^T\Sigma^{-1}\bZ\right)^{m/2}\over \left(1+\bZ^T\Sigma^{-1}\bZ
b^2/v\right)^{(v+m)/2}v^{m/2}\Gamma(m/2)\Gamma(v/2)}\right].
\end{eqnarray*}

%\newpage

%\section{Online Appendix}

\end{document}